\documentclass[10pt,a4paper, reqno]{amsart}
\usepackage{mathrsfs}

\usepackage{amsmath,amsfonts,verbatim}
\usepackage{latexsym}
\usepackage{amssymb,leftidx}
\usepackage{extarrows}
\usepackage{overpic}
\usepackage{color}
\usepackage{epsfig}
\usepackage{subfigure}
\usepackage{tikz}
\usepackage{bbm}
\usepackage{tikz, caption}
\usepackage[font=small,labelfont=bf]{caption}

\usepackage[colorlinks=true,backref=page]{hyperref}
\hypersetup{urlcolor=blue, citecolor=blue}

\usepackage{amssymb}
\usepackage{mathrsfs}
\usepackage{amscd}
\usepackage{cite}

\newcommand\R{\mathbb{R}}

\newcommand\Z{\mathbb{Z}}
\newcommand\N{\mathbb{N}}

\newcommand\A{{\bf A}}

\newcommand\LL{\mathcal{L}}

\usepackage{graphicx}
\usepackage{float}

\numberwithin{equation}{section}
\newtheorem{proposition}{Proposition}[section]

\newtheorem{lemma}{Lemma}[section]
\newtheorem{theorem}{Theorem}[section]

\newtheorem{remark}{Remark}[section]

\begin{document}
\title[$L^p$-estimates for magnetic wave]{$L^p$-estimates for the 2D wave equation\\ in the scaling-critical magnetic field}

\author{Jialu Wang}
\address{Department of Mathematics, Beijing Institute of Technology, Beijing, China, 100081;} \email{jialu\_wang@bit.edu.cn}

\author{Fang Zhang}
\address{Department of Mathematics, Beijing Institute of Technology, Beijing, China, 100081;} \email{zhangfang@bit.edu.cn}

\author{Junyong Zhang}
\address{Department of Mathematics, Beijing Institute of Technology, Beijing, China, 100081;} \email{zhang\_junyong@bit.edu.cn}

\author{Jiqiang Zheng}
\address{Institute of Applied Physics and Computational Mathematics, Beijing 100088}
\email{zheng\_jiqiang@iapcm.ac.cn}

\begin{abstract}
In this paper, we study the $L^{p}$-estimates for the solution to the $2\mathrm{D}$-wave equation with a scaling-critical magnetic potential. Inspired by the work of \cite{FZZ}, we show that the operators $(I+\mathcal{L}_{\mathbf{A}})^{-\gamma}e^{it\sqrt{\mathcal{L}_{\mathbf{A}}}}$ is bounded in $L^{p}(\mathbb{R}^{2})$ for $1<p<+\infty$ when $\gamma>|1/p-1/2|$ and $t>0$, where $\mathcal{L}_{\mathbf{A}}$ is a magnetic Schr\"odinger operator. In particular, we derive the $L^{p}$-bounds for the sine wave propagator $\sin(t\sqrt{\mathcal{L}_{\mathbf{A}}})\mathcal{L}^{-\frac12}_{\mathbf{A}}$.
The key ingredients are the construction of the kernel function and the proof of the pointwise estimate for an analytic operator family $f_{w,t}(\mathcal{L}_{\mathbf{A}})$.

\end{abstract}
 \maketitle

\begin{center}
\begin{minipage}{120mm}
   { \small {{\bf Key Words:} $L^{p}$-estimates; Scaling-critical magnetic field; Aharonov-Bohm potential; Wave equation.}
   }\\
\end{minipage}
\end{center}
\section{introduction}

In this paper, we investigate the $L^{p}$-estimates of the solution for the wave equations with a scaling-critical magnetic potential.
More precisely, we consider the following wave equation
\begin{align}\label{equation}
\begin{cases}
\partial_{tt}u+\mathcal{L}_{\mathbf{A}}u=0,\ (t, x)\in \mathbb{R}\times \mathbb{R}^{2},\\
u(0, x)=f(x),\ \partial_{t}u(0, x)=g(x),
\end{cases}
\end{align}
where the magnetic singular Schr\"odinger operator is given by
\begin{align}\label{LA}
\mathcal{L}_{\mathbf{A}}=\Big(i\nabla+\frac{\mathbf{A}(\hat{x})}{|x|}\Big)^{2},\quad x\in \mathbb{R}^{2}\backslash\{0\}.
\end{align}
Here $\hat{x}=\frac{x}{|x|}\in \mathbb{S}^1$ and $\mathbf{A}\in W^{1,\infty}(\mathbb{S}^1;\R^2)$ satisfies the transversality condition
\begin{equation}\label{transcondition}
\mathbf{A}(\hat{x})\cdot \hat{x}=0,\quad \text{for all} \quad \hat{x}\in\mathbb{S}^1.
\end{equation}
One typical example is the following {\it Aharonov-Bohm} potential
\begin{equation}\label{AB}
\mathbf{A}(\hat{x})=\alpha\left(\frac{-x_{2}}{|x|}, \frac{x_{1}}{|x|}\right),\ \alpha\in \mathbb{R},
\end{equation}
which was introduced initially in \cite{AB}, in the context of Schr\"odinger dynamics, to show that scattering effects can even occur in regions in which the electromagnetic field is absent (see also \cite{PT}). Fanelli and the last two authors \cite{FZZ} have established the Strichartz estimates by constructing the propagator $\frac{\sin(t\sqrt{\mathcal{L}_{\mathbf{A}}})}{\sqrt{\mathcal{L}_{\mathbf{A}}}}$ (based on {\it Lispschitz-Hankel integral formula}). In this paper, we aim to study the $L^p$ estimates for the wave propagator perturbed by the scaling-critical magnetic singular potentials.\vspace{0.2cm}

The study of $L^p$-estimates for the wave equation traces back to Euclidean space $\R^n$, classically going back to the work of \cite{M, P}. Specifically, Peral \cite{P} and Miyachi \cite{M} established the sharp range of $p$, namely $|\frac{1}{p}-\frac{1}{2}|\leq \frac{1}{n-1}$, for which $\frac{\sin(t\sqrt{-\Delta})}{\sqrt{-\Delta}}$ is bounded on $L^p(\R^n)$. The $L^p$-bounds for the operators $(Id-\Delta)^{-\gamma}e^{it\sqrt{-\Delta}}$ for $\gamma\geq (n-1)|1/p-1/2|$ are discussed in \cite{M} and the $L^p\rightarrow L^q$ estimates were explored by Strichartz \cite{S}.
It is important to note that these results are proven for the operator $-\Delta_{\R^n}$, where classical Fourier analysis techniques are applicable. However, analogous problems for variable coefficient Schrödinger operators, such as those perturbed by potentials, introduce additional complexities.
For instance, the Hermite operator $-\Delta+|x|^2$ was investigated by \cite{NT}, and more general operators of the form $-\Delta+V$ were studied by Zhong \cite{Z}.  K. Jotsaroop and S. Thangavelu \cite{KT} later demonstrated that the operator $\frac{\sin\sqrt{G}}{G}$, in which the Grushin operator $G=-\Delta-|x|^2\partial_{t}^2$, is bounded on $L^p(\R^{n+1})$ for every $p$ satisfying $|\frac{1}{p}-\frac{1}{2}|<\frac{1}{n+2}$, thereby achieving the $L^p$ estimate for the solution. Correspondingly, they proved $L^p\to L^2$ estimate for wave equation with the Grushin operator via employing the $L^p$ boundedness of Riesz transform for the Grushin operator in subsequent article \cite{TN}. The $L^p$-boundedness of oscillating multipliers, on certain wide classes of rank one locally symmetric spaces, was deduced in \cite{E}.  Li and Lohoue in \cite{LN, L} proved the $L^p$-estimates for the wave equation on cone $C(N)$.  We also refer to $L^p$-estimates for the wave equation on other settings:
Heisenberg groups, as discussed in \cite{MS2};
symmetric spaces of non-compact type with real rank one, in \cite{GM,Ionescu}; Damek-Ricci spaces, covered in \cite{MV}; variable coefficient Fourier integral operators, \cite{HPR, LRSY,HR, R}.
For extra results on the wave equation in manifolds with conical singularities, refer to \cite{Mewun,MS1,SS1,SS2}. In this paper, we focus on the wave equation with a scaling-critical singular magnetic potential.\vspace{0.2cm}

Regarding the operator $\mathcal{L}_{\mathbf{A}}$, there are several related results worth mentioning.
In \cite{FFFP}, Fanelli, etc. proved the time-decay estimate
\begin{equation}\label{time-decay}
\|e^{it\mathcal{L}_{\mathbf{A},a}}\|_{L^1(\R^2)\rightarrow L^{\infty}(\R^2)}\leq|t|^{-1}
\end{equation}
for Schr\"odinger equation with scaling invariant electromagnetic potential in which
\begin{equation}
\mathcal{L}_{\mathbf{A},a}=\mathcal{L}_{\mathbf{A}}+\frac{a(\hat{x})}{|x|^2},\quad \text{with}\quad a\in W^{1,\infty}(\mathbb{S}^1,\R).
\end{equation}
So the Strichartz estimates can be obtained directly by the Keel-Tao argument. However, the argument in \cite{FFFP} breaks down for wave and Klein-Gordon equations due to the lack of  pseudoconformal invariance (which was used for Schr\"odinger equation). Very recently, Fanelli, and the last two authors \cite{FZZ}
proved the Strichartz estimate for wave equation by constructing the propagator $\sin(t\sqrt{\mathcal L_{{\A},0}})/\sqrt{\mathcal L_{{\A},0}}$ (based on \emph{Lipschitz-Hankel integral formula}) and showing the local smoothing estimates. Gao, Yin, and the last two authors \cite{GYZZ} constructed the spectral measure and further proved the time decay and Strichartz estimates of Klein-Gordon equation. Based on the spectral measure, we studied the $L^p-L^q$- resolvent estimates for the Schr\"odinger operator $\mathcal{L}_{\mathbf{A}}$ in \cite{FZZ1} and proved the uniform weighted $L^2$-estimates  in \cite{GWZZ}. \vspace{0.2cm}

Recently, we start a new program to study the $L^p$-bounds for spectral multipliers associated with the self-adjoint (Friedrichs extension) operator $\mathcal{L}_{\A}$. More precisely, using the spectral theorem, the ``abstract spectral multipliers'' $F(\mathcal{L}_{\A})$ can be defined for any bounded Borel function $F$, and they act continuously on $L^2$. Our problem is to find some necessary conditions on the function $F$ to ensure that the operator $F(\mathcal{L}_{\A})$ extends as a bounded operator for some range of $L^p$ spaces for $p\neq2$. The similar problem associated with various self-adjoint operators has attracted considerable attention over the past forty years, now we focus on the magnetic Schr\"odinger operator $\LL_{\A}$. Additionally, certain families of functions $F$ have been examined within the framework of spectral multiplier theory. Notably, the bounds for the Littlewood-Paley operator and its multipliers can be derived from these abstract results, particularly given the Gaussian upper bounds for the heat kernel of $\mathcal{L}_{\A}$. For further details, we refer the reader to \cite{FZZ}. In \cite{MYZ}, Miao, Yan and the third author proved the sharp $L^p$-boundedness of the Bochner-Riesz means associated with the operator $\mathcal{L}_{\mathbf{A}}$. In that paper, we find the magnetic Aharonov-Bohm effect play roles. For technique issue, the incident wave is different from the usual plane wave $e^{ix\cdot\xi}$, which leads to that the modified factor $e^{\pm i\alpha(\theta_1-\theta_2)}$ appears in the kernel of Bochner-Riesz means because of the long-range property of the potential $\frac{{\A}(\hat{x})}{|x|}$. In addition, the Heaviside step function $\mathbbm{1}_{I}(\theta_1-\theta_2)$ appears in the kernel of Bochner-Riesz means,
which is another obstacle to efficiently exploit the oscillation behavior of the kernel.\vspace{0.2cm}


In this paper,  we are interested in proving the $L^p$-bounds for oscillatory integrals of the form
\begin{equation}\label{osc-oper}
e^{i\mathcal{L}_{\A}^{\frac m2}}(Id+\mathcal{L}_{\A}^{\frac m2})^{-\gamma},
\end{equation}
and we intend to determine the assumption on the optimal regularity index $\gamma$ and $p$ such that the operator \eqref{osc-oper} is $L^p$-bounded for either $m=2$ or $m=1$, which are corresponding to Schr\"odinger and wave equations respectively. For the Schr\"odinger, i.e. $m=2$, the \cite[Theorem 1]{CDLY} and \cite[Theorem 7.19]{Ou}  show that the operator
\eqref{osc-oper} is $L^p$-bounded for $\gamma\geq 2|1/2-1/p|$ (This result is sharp including the endpoint.) if one could prove
the Gaussian upper bounds of the heat kernel
\begin{align}\label{equ:pcontr}
 |e^{-t\LL_{{\A}}}(x,y)|\lesssim \,& t^{-1}e^{-\frac{|x-y|^2}{4t}},\quad \forall\;t>0,
\end{align}
which was prove in \cite[Proposition 3.1, Proposition 3.2]{FZZ} by constructing the heat kernel for the Schr\"odinger operator $\mathcal{L}_{\mathbf{A}}$.
We also refer to \cite{Ko, MOR} for heat kernel with Aharonov-Bohm potentials \eqref{AB} and reference therein.
For the wave equation, one can use \cite[Theorem 7.19]{Ou} and \eqref{equ:pcontr} to show the operator \eqref{osc-oper} with $m=1$ is $L^p$-bounded for $\gamma>2|1/2-1/p|$ but is far from the sharp result.
Motivated by this observation, in this paper, we aim to improve the regularity assumption to $\gamma>|1/2-1/p|$ which is sharp up to the endpoint $\gamma=|1/2-1/p|$.
More precisely, we prove

\begin{theorem}\label{Lpestimate}
Let $\LL_{\A}$ be given magnetic Schr\"odinger operator in \eqref{LA} and let $\gamma>0$ and $1< p<+\infty$ satisfy $|\frac{1}{p}-\frac{1}{2}|<\gamma$. Then, there exists a constant $C(p, \gamma)>0$ such that for $t>0$ and for all $f\in L^{p}(\mathbb{R}^{2})$
\begin{equation}
\|(1+\mathcal{L}_{\mathbf{A}})^{-\frac{\gamma}{2}}e^{it\sqrt{\mathcal{L}_{\mathbf{A}}}}f\|_{L^{p}(\R^2)}\leq C(p, \gamma)(1+t)^{\gamma}\|f\|_{L^{p}(\R^2)}.
\end{equation}
\end{theorem}

\begin{remark}
This result is sharp up to the endpoint $\gamma=|\frac1p-\frac12|$. In particular, $\gamma=|\frac1p-\frac12|$,
the argument fails due to the fact that there is no room allowing one choose $\epsilon<1$ in \eqref{est:fw} below.
\end{remark}

It has been known that the solution to \eqref{equation} can be represented by
\begin{equation}
u(t,\cdot)=\cos(t\sqrt{\mathcal{L}_{\mathbf{A}}})f(\cdot)+\frac{\sin(t\sqrt{\mathcal{L}_{\mathbf{A}}})}{\sqrt{\mathcal{L}_{\mathbf{A}}}}g(\cdot).
\end{equation}
Our result related to the wave equation \eqref{equation} is concerned with the following.
\begin{theorem}\label{sinLA}
Let $\LL_{\A}$ be given magnetic Schr\"odinger operator in \eqref{LA}. There exists a constant $C>0$ such that, for any $1\leq p\leq+\infty,$
\begin{equation}\label{est:sinLA}
 \left\|\frac{\sin(t\sqrt{\mathcal{L}_{\mathbf{A}}})}{\sqrt{\mathcal{L}_{\mathbf{A}}}}g\right\|_{L^{p}(\mathbb{R}^{2})}\leq C |t|\cdot\Vert g\Vert_{L^{p}(\mathbb{R}^{2})},\quad \forall\; t\neq 0.
\end{equation}
\end{theorem}

Now we sketch our proof here. Our strategy is to combine the method in \cite{L} and the construction of the kernels in \cite{FZZ, GYZZ}.
As \cite{L, MS1} did, we  approximate the half-wave multiplier $e^{it|\xi|}$ by $\cos(t|\xi|)$ which is closed to the Bessel function as $\rho=|\xi|\to +\infty$ in view of the asymptotics \cite{W}
\begin{equation}
\Gamma(\nu+\frac12)\frac{J_{\nu-\frac12}(\rho)}{\rho^{\nu-\frac12}}=\begin{cases}1+o(1),\quad \rho\to 0,\\
(2/\pi)^{1/2} \rho^{-\nu}\cos(\rho-\frac{\pi\nu}2), \quad \rho\to \infty.\end{cases}
\end{equation}
Hence we consider the family of analytic operators defined by
\begin{equation}\label{def:ana-oper}
f_{w,t}(\mathcal{L}_{\mathbf{A}})=\big(\frac{\pi}{2}\big)^{\frac{1}{2}}(t\sqrt{\mathcal{L}_{\mathbf{A}}})^{w-1}J_{1-w}(t\sqrt{\mathcal{L}_{\mathbf{A}}})\ ,
\end{equation}
where $w=\epsilon+iy(\frac{3}{2}-\epsilon)$ with $y\in \mathbb{R}$ and $\frac{1}{2}<\epsilon<1.$
The first key ingredient in the proof is the explicit formulas \eqref{formula:Ma} below for the integral expressions
\begin{align}
\int_{0}^{+\infty} t^{1-\mu} J_{\mu}(at) J_{\lambda}(bt) J_{\lambda}(ct) \, dt,
\end{align}
which are due to Macdonald \cite{Ma}. The second ingredient is the summation in the spectrums to construct the kernels, as did in \cite{FZZ, GYZZ}.
Hence we can prove the $L^p$-bounds for the analytic operators $f_{w,t}(\mathcal{L}_{\mathbf{A}})$.
\begin{theorem}\label{operator-family}
Let $\frac{1}{2}<\epsilon<1$ and $ 1\leq p\leq+\infty$. Then, there exists a constant $C(\epsilon,p)>0$ such that
\begin{equation}\label{est:fw}
\|f_{w,t}(\mathcal{L}_{\mathbf{A}})g\|_{L^p(\R^2)}\leq C(\epsilon,p)\frac{e^{-y\pi}}{|\Gamma(1-w)|}\| g\|_{L^p(\R^2)}, \quad\forall\; g\in L^{p}(\mathbb{R}^{2}),
\end{equation}
where $w=\epsilon+iy$ with $y\in \mathbb{R}$.
\end{theorem}
Eventually, we can prove Theorem \ref{Lpestimate} as a consequence of this theorem.\vspace{0.2cm}

It is important to note the differences between \cite{L, LN} and the current work, even though we are inspired by them.  In \cite{L}, Li employed the parametrix method to demonstrate the boundedness of the corresponding family of analytic operators $f_{w,t}(\Delta_{g})$ with the dimension $n\geq3$, in which $\Delta_{g}$ is the Laplace-Beltrami operator on the cone manifolds $C(N)$. However, this approach is not applicable in our case which involves a two-dimensional setting.
As a result, rather than employing their parametrix method, we precisely determine the kernel function of the analytic operator family $f_{w,t}(\mathcal{L}_{\mathbf{A}})$ via implementing functional calculus and fulfill the pointwise estimate. Moreover, the authors of \cite{L} and \cite{LN} used Stein's interpolation theorem of an analytical operator family defined by Bessel functions to justify Theorem \ref{sinLA} for $\Delta_{g}$ on the cone manifolds $C(N)$ with the dimension $n\geq3$. Nevertheless, it isn't suitable for the magnetic operator $\mathcal{L}_{\mathbf{A}}$ due to the range limit of $w$. Thus, we instead explore the kernel estimates to prove Theorem \ref{sinLA}.\vspace{0.2cm}

The structure of this paper is organized as the following. In Section \ref{analysis}, we provide some analysis results associated with the operator $\mathcal{L}_{\mathbf{A}}$. In Section \ref{proof-family}, Theorem \ref{operator-family} is demonstrated by showing pointwise estimate for kernel of analytic operator family $f_{w,t}(\mathcal{L}_{\mathbf{A}})$. In Section \ref{proof-theorem} and Section \ref{proof-sinLA}, we prove Theorem \ref{Lpestimate} and Theorem \ref{sinLA} respectively. Finally, the classical property of the Bessel function is prepared in the Appendix section.\vspace{0.2cm}

{\bf Acknowledgments:}\quad J. Zhang was supported by National Natural Science Foundation of China (12171031) and Beijing Natural Science Foundation (1242011). J. Zheng was supported by National key R\&D program of China (Nos. 2021YFA1002500, 2020YFA0712900), Beijing Natural Science Foundation (No. 1222019), NSFC (Nos. 12271051, 11831004).

\section{Preliminaries}\label{analysis}
In this preliminary section, we firstly introduce some notations used in this article. In addition, we recall and provide some analysis tools, including the spectral property and spectral multiplier theorem associated with the operator $\mathcal{L}_{\mathbf{A}}$.
\subsection{Notation} Throughout this article, we utilize $C$ to denote the universal constant. We say that $A\lesssim B$ if $A\leq CB$. Additionally, we use the notation $C(B)$ to refer to a constant that depends on $B$. For convenience, we write $L^p(\R^2)\rightarrow L^p(\R^2)$ as $p\rightarrow p$. Let $\Omega$ be a nonempty subset of $\R^n$, we will denote $C^{m}(\Omega)$ as the space of continuous functions defined on $\Omega$, while $B(\Omega)$ represents the set of bounded functions on $\Omega$. Furthermore, we employ the notation $C^{m,\lambda}(\Omega)$ to denote the usual H\"older space, in which $m$ is a nonnegative integer and $\lambda\in(0,1]$.

The H\"older space is defined by
\begin{equation}
C^{m,\lambda}(\Omega)=\{u\in C_{b}^{m}(\Omega): \|u\|_{C^{m,\lambda}(\Omega)}<\infty\},
\end{equation}
where
\begin{equation}
\|u\|_{C^{m,\lambda}(\Omega)}=\|u\|_{C^{m}(\Omega)}+\max_{|\alpha|=m}[\partial^{\alpha}u]_{\lambda,\Omega},
\end{equation}
\begin{equation}
C_{b}^{m}(\Omega)=\{u\in C^{m}(\Omega): \partial^{\alpha}u\in B(\Omega), \forall |\alpha|\leq m\},
\end{equation}
and
\begin{equation}
[u]_{\lambda,\Omega}=\sup_{x,y\in\Omega, x\neq y}\frac{|u(x)-u(y)|}{|x-y|^{\lambda}}.
\end{equation}
By the definition of H\"older space, for $s>0$, we have
\begin{equation}
C^{s}(\Omega)=C^{[s], s-[s]}(\Omega)
\end{equation}
in which $[s]$ denotes the integer part of $s$.

Also, we define the Hankel transform of order $\nu$ as follows
\begin{equation*}
(\mathcal{H}_{\nu}f)(\rho,\theta)=\int_{0}^{\infty}J_{\nu}(r\rho)f(r, \theta)rdr,
\end{equation*}
where the Bessel function of order $\nu$ is given by
\begin{equation*}
J_{\nu}(r)= \frac{(r/2)^{\nu}}{\Gamma(\nu+\frac{1}{2})\Gamma(1/2)}\int_{-1}^{1}e^{isr}(1-s^{2})^{(2\nu-1)/2}ds, \quad \nu>-1/2, \ r>0.
\end{equation*}

\subsection{Functional calculus associated with the operator $\mathcal{L}_{\mathbf{A}}$}
In this subsection, inspired by Cheeger-Taylor \cite{CT}, we study the harmonic analytical features of the operator $\mathcal{L}_{\mathbf{A}}$ and recall the functional calculus connected with the operator $\mathcal{L}_{\mathbf{A}}$.

Recall the operator $\mathcal{L}_{\mathbf{A}}$ in \eqref{LA} given by
\begin{align}\label{def:LA}
\mathcal{L}_{\mathbf{A}}=-\Delta+|x|^{-2}(|\mathbf{A}(\hat{x})|^2+i \text{div}_{\mathbb{S}^1}\mathbf{A}(\hat{x}))+2i|x|^{-1}\mathbf{A}(\hat{x})\cdot \nabla.
\end{align}
Using \eqref{transcondition} and polar coordinates, we express $\mathcal{L}_{\mathbf{A}}$ as follows
\begin{equation}\label{def:LbfA} \mathcal{L}_{\mathbf{A}}=-\partial_{r}^{2}-\frac{1}{r}\partial_{r}+\frac{L_{\mathbf{A}}}{r^{2}},  \end{equation}
in which the operator $L_{\mathbf{A}}$ is given by
\begin{align}
L_{\mathbf{A}}&=(i\nabla_{\mathbb{S}^{1}}+\mathbf{A}(\hat{x}))^{2},\quad \hat{x}\in \mathbb{S}^{1}\nonumber\\
&=-\Delta_{\mathbb{S}^{1}}+(|\mathbf{A}(\hat{x})|^{2}+i\mathrm{d}\mathrm{i}\mathrm{v}_{\mathbb{S}^{1}}\mathbf{A}(\hat{x}))+2i\mathbf{A}(\hat{x})\cdot\nabla_{\mathbb{S}^{1}}.
\end{align}
Let $\hat{x}=(\cos\theta,\sin\theta)$, then we have
$$
\partial_{\theta}=-\hat{x}_{2}\partial_{\hat{x}_{1}}+\hat{x}_{1}\partial_{\hat{x}_{2}},\quad \partial_{\theta}^{2}=\Delta_{\mathbb{S}^{1}}.
$$
We define $\alpha(\theta) : [0,2\pi]\rightarrow \mathbb{R}$ such that
\begin{equation}
\alpha(\theta)=\mathbf{A}(\cos\theta, \sin\theta)\cdot(-\sin\theta, \cos\theta).
\end{equation}
By applying the transversality condition \eqref{transcondition}, we can obtain
\begin{equation}
\mathbf{A}(\cos\theta, \sin\theta)=\alpha(\theta)(-\sin\theta, \cos\theta),\quad \theta\in[0,2\pi].
\end{equation}
Consequently, we further rewrite $L_{\mathbf{A}}$ as
\begin{align}
L_{\mathbf{A}}
=(i\partial_{\theta}+\alpha(\theta))^{2}.
\end{align}

Next, we define the magnetic flux along the sphere
\begin{equation}\label{alphadefine}
\alpha:=\Phi_{\mathbf{A}}=\frac{1}{2\pi}\int_{0}^{2\pi}\alpha(\theta)d\theta.
\end{equation}
According to the result of Laptev and Weidl \cite{LW}, the operator $i\partial_{\theta}+\alpha(\theta)$ with domain $H^{1}(\mathbb{S}^{1})$ in $L^{2}(\mathbb{S}^{1})$ possesses eigenvalues $\nu_{k}=k+\Phi_{\mathbf{A}}, k\in \mathbb{Z}$, and the corresponding eigenfunctions are given by
\begin{equation}\label{eig-LA} \varphi_{k}(\theta)=\frac{1}{\sqrt{2\pi}}e^{-i(\theta(k+\alpha)-\int_{0}^{\theta}\alpha(\theta')d\theta')}.
\end{equation}
Thus, we have the relation
$$
L_{\mathbf{A}}\varphi_{k}(\theta)=(k+\alpha)^{2}\varphi_{k}(\theta).
$$
Based on the eigenfunctions, we can acquire the orthogonal decomposition $$L^2(\mathbb{S}^{1})=\bigoplus_{k\in\N}\mathcal{H}^k,\quad \mathcal{H}^k=\text{span}\{\varphi_{k}(\theta)\}.$$
Then, given \eqref{def:LbfA}, we can conclude that, on each space $\mathcal{H}^{k}=\mathrm{s}\mathrm{p}\mathrm{a}\mathrm{n}\{\varphi_{k}\}$, the action of the operator is expressed as
\begin{equation}
\mathcal{L}_{\mathbf{A}}=-\partial_{r}^{2}-\frac{1}{r}\partial_{r}+\frac{(k+\alpha)^{2}}{r^{2}}.
\end{equation}
with eigenvalues $\nu=\nu_{k}=|k+\alpha|$ and $k\in \mathbb{Z}$.

Also, for any $f\in L^2(\R^2)$, we can express $f$ in the form of separating variables
\begin{equation}\label{sep.v}
f(x)=\sum_{k\in\N} a_{k}(r)\varphi_k(\theta),
\end{equation}
where
\begin{equation*}
 a_{k}(r)=\int_{0}^{2\pi}f(r,\theta)
\overline{\varphi_k(\theta)}\, d\theta.
\end{equation*}

%
%
%
%

Ultimately, we shortly recall the functional calculus for well-behaved functions $F$ (see \cite{Taylor}),
\begin{equation*}
F(\mathcal{L}_{\mathbf{A}})f(r_{1},\theta_{1})=\sum_{k\in \mathrm{Z}}\varphi_{k}(\theta_{1})\int_{0}^{\infty}F(\rho^{2})J_{\nu_{k}}(r_{1}\rho)b_{k}(\rho)\rho d\rho,
\end{equation*}
where $b_{k}(\rho)=(\mathcal{H}_{\nu_{k}}a_{k})(\rho)$ and $f(r_{1},\theta_{1})=\sum_{k\in \mathbb{Z}}a_{k}(r_{1})\varphi_{k}(\theta_{1})$. Thus,
\begin{align}\label{FKmuk}
F(\mathcal{L}_{\mathbf{A}})f(r_{1},\theta_{1})&=\sum_{k\in \mathbb{Z}}\varphi_{k}(\theta_{1})\int_{0}^{\infty}F(\rho^{2})J_{\nu_{k}}(r_{1}\rho)b_{k}(\rho)\rho d\rho\nonumber\\
&=\int_{0}^{\infty}\int_{0}^{2\pi}f(r_{2},\theta_{2})K(r_{1}, \theta_{1}, r_{2}, \theta_{2})r_{2}dr_{2}d\theta_{2},
\end{align}
in which we have used the following
$$
K(r_{1},\theta_{1}, r_{2}, \theta_{2})=\sum_{k\in \mathbb{Z}}\varphi_{k}(\theta_{1})\overline{\varphi_{k}(\theta_{2})}K_{\nu_{k}}(r_{1}, r_{2}) ,
$$
and
\begin{equation*}
K_{\nu_{k}}(r_{1},r_{2})=\int_{0}^{\infty}F(\rho^{2})J_{\nu_{k}}(r_{1}\rho)J_{\nu_{k}}(r_{2}\rho)\rho d\rho.
\end{equation*}
Now we take $F=f_{w,t}$ in \eqref{def:ana-oper}, so we need to study the kernel
\begin{equation}\label{Kmuk}
K_{\nu_{k}}(r_{1},r_{2})=\Big(\frac{\pi}{2}\Big)^{\frac{1}{2}}\int_{0}^{\infty}(t\rho)^{w-1}J_{1-w}(t\rho) J_{\nu_{k}}(r_{1}\rho)J_{\nu_{k}}(r_{2}\rho)\rho d\rho.
\end{equation}

To this end, we adopt some notations from \cite{W} for special functions. We denote the Gamma function as $\Gamma(z)$. Also, we use $P_{\lambda}^{\mu}(x)$ and $Q_{\lambda}^{\mu}(z)$ to represent the Legendre functions when $|x|<1$ and $|z|>1$ respectively. Now, we will present the following facts about these special functions which were used in \cite{L}:

\underline{\textbf{Fact I:}} Let $\mathrm{Re}\mu\neq-\frac{1}{2}, \mathrm{Re} \lambda\neq-\frac{1}{2}$ and $a, b, c>0$. Define
$$A= \arccos\Big(\frac{b^{2}+c^{2}-a^{2}}{2bc}\Big),\quad \mathcal{A}=\cosh^{-1}\Big(\frac{a^{2}-b^{2}-c^{2}}{2bc}\Big).$$
Then, by \cite[p.412]{W}, we have
\begin{align}\label{formula:Ma}
&\int_{0}^{+\infty} t^{1-\mu} J_{\mu}(at) J_{\lambda}(bt) J_{\lambda}(ct) \, dt \nonumber\\
&=
\begin{cases}
0, & \text{if } a < |b - c|, \\
\frac{(bc)^{\mu - 1} \sin^{\mu - \frac{1}{2}} A}{(2\pi) 2 a^{\mu}} P_{\lambda - \frac{1}{2}}^{\frac{1}{2} - \mu}(\cos A), & \text{if } |b - c| < a < b + c, \\
\frac{(bc)^{\mu - 1} \sinh^{\mu - \frac{1}{2}} \mathcal{A}}{\left( \frac{\pi^{3}}{2} \right)^{\frac{1}{2}} a^{\mu}} \cos(\pi \lambda) Q_{\lambda - \frac{1}{2}}^{\frac{1}{2} - \mu}(\cosh \mathcal{A}), & \text{if } a > b + c.
\end{cases}
\end{align}

\underline{\textbf{Fact II:}} Let $ 0<\theta<\pi$ and $\mathrm{Re} \mu_{*}<\frac{1}{2}$. Then, by \cite[p.188]{M}, we can see
\begin{equation}\label{gamma}
\begin{split}
 \Gamma(\frac{1}{2}-\mu_{*})&P_{\lambda}^{\mu_{*}}(\cos\theta)\\
&=\Big(\frac{\pi}{2}\Big)^{-\frac{1}{2}}(\sin\theta)^{\mu_{*}}\int_{0}^{\theta}(\cos s-\cos\theta)^{-\mu_{*}-\frac{1}{2}}\cos[(\lambda+\frac{1}{2})s]\,ds.
\end{split}
\end{equation}

\underline{\textbf{Fact III:}} Let $\mathrm{Re}(\lambda+\mu_{*}+1)>0$ and $\mathrm{Re}\mu_{*}<\frac{1}{2}$. Then, one has
\begin{align}\label{gammaQ}
\Gamma(\frac{1}{2}-\mu_{*})&Q_{\lambda}^{\mu_{*}}(\cosh\theta)\nonumber
\\&=\Big(\frac{\pi}{2}\Big)^{\frac{1}{2}}e^{i\pi\mu_{*}}(\sinh\theta)^{\mu_{*}}\int_{\theta}^{+\infty}(\cosh s-\cosh\theta)^{-\mu_{*}-\frac{1}{2}}e^{-(\lambda+\frac{1}{2})s}\,ds.
\end{align}
See \cite[p.186]{M}.

Therefore, combining \eqref{formula:Ma}, \eqref{gamma}, \eqref{gammaQ} and the construction in \cite{FZZ, GYZZ}, we will construct the kernel of $f_{w,t}(\LL_{\A})$ in Section \ref{proof-family}.

\subsection{Some results about spectral multipliers}In this subsection, we present some results concerning the spectral multipliers derived from the heat kernel estimate
\begin{equation}\label{estimate:heatkernel}
\left|e^{-t \mathcal{L}_{\mathbf{A}}}(x, y)\right| \lesssim \frac{1}{t} \exp \left(-\frac{|x-y|^2}{4 t}\right), \quad \forall t>0,
\end{equation}
which was proved in \cite{FZZ}.

\begin{lemma}\cite[Theorem 7.23]{Ouhabaz}\label{theorem-spectral}
Assume that $(X, \rho, \mu)$ satisfies the doubling property. Let $\Omega$ be an open subset of $X$ and assume that $A$ is a non-negative self-adjoint operator on $L^2(\Omega,\mu)$ with a heat kernel $p(t,x,y)$ satisfying the Gaussian upper bound. Denote by $\varphi$ a non-negative $C_{c}^{\infty}$ function satisfying
\begin{equation}
\textrm{supp}(\varphi)\subseteq[\frac{1}{4},1] \quad\text{and} \quad \sum_{n\in\mathbb{Z}}\varphi(2^{-n}\lambda)=1, \quad\text{for all} \ \lambda>0.
\end{equation}
Let $f: [0,\infty)\rightarrow \mathbb{C}$ be a bounded function such
that
\begin{equation}\label{estimate-condition1}
\sup_{t>0}\|\varphi(\cdot)f(t\cdot)\|_{C^s}<\infty
\end{equation}
for some $s>\frac{d}{2}$. Then $f(A)$ is of weak type $(1,1)$ and is bounded on $L^p(\Omega,\mu)$ for all $p\in(1,\infty)$. In addition,
\begin{equation}
\|f(A)\|_{L^1(\Omega,\mu)\rightarrow L^{1,w}(\Omega,\mu)}\leq C_{s}\big(\sup_{t>0}\|\varphi(\cdot)f(t\cdot)\|_{C^s}+|f(0)|\big)
\end{equation}
for some positive constant $C_{s}$, independent of $f$.
\end{lemma}

\begin{lemma}\label{lem:multiplier}
Let $\LL_{\A}$ be the magnetic Schr\"odinger operator in \eqref{LA}.

$\bullet$ $L^p$-bounds for imaginary powers operator.  For all $y\in\R$, then the imaginary powers $(\LL_{\A})^{iy}$ satisfy the $(1,1)$ weak type estimate
\begin{equation}\label{im-bounds}
\big\|(\LL_{\A})^{iy}\big\|_{L^1(\R^2)\to L^{1,\infty}(\R^2)}\leq C(1+|y|)^{n/2}.
\end{equation}
In particular,  there exists a constant $C>0$ such that
\begin{equation}\label{spectral-LA}
\|(\LL_{\A})^i\|_{L^p(\R^2)\to L^{p}(\R^2)} \leq C, \quad\forall 1<p<+\infty.
\end{equation}

$\bullet$ Mikhlin-H\"ormander multiplier estimates. Assume that $m:[0,+\infty) \rightarrow \mathbb{R}$ satisfies
\begin{equation}
\sum_{j=0}^2 \sup _{s \geq 0} s^j|m^{(j)}(s)|<+\infty,
\end{equation}
then, for all $1<p<+\infty$ and $t>0$,
\begin{equation}\label{estimate:spectraltLA}
\|m(t \sqrt{\LL_{\A}})\|_{L^p(\R^2)\to L^{p}(\R^2)} \leq C \sum_{j=0}^2 \sup _{s \geq 0} s^j|m^{(j)}(s)| .
\end{equation}

$\bullet$ Let $\sigma>1$ and $m:[0,+\infty) \to \mathbb{R}$ obeys
\begin{equation}\label{estimate-condition}
|m^{(j)}(s)| \leq C (1+s)^{-\sigma}, \quad\forall s \geq 0 \text { and } 0 \leq j \leq 2,
\end{equation}
then there exists a constant $C(\sigma)>0$ such that
\begin{equation}\label{ResultIII}
\|m(t \sqrt{\LL_{\A}})\|_{L^p(\R^2)\to L^{p}(\R^2)} \leq C(\sigma), \quad \forall 1<p<+\infty, \ t>0 .
\end{equation}

\end{lemma}

\begin{proof} The proof is based on the heat kernel estimate \eqref{estimate:heatkernel}. The inequality \eqref{im-bounds} follows from \cite[Theorem 2]{SW}.
Hence one can obtain \eqref{spectral-LA} by interpolation, see \cite[Corollary 7.24]{Ouhabaz}. The proof of \eqref{estimate:spectraltLA} is standard if one has
the Gaussian upper bounds  \eqref{estimate:heatkernel}, we refer to \cite[(7.69)]{Ouhabaz}.

We now turn to the proof of \eqref{ResultIII}. For $\sigma\geq2$, \eqref{ResultIII} follows directly from \eqref{estimate:spectraltLA}. Therefore, we will focus on the case where $1<\sigma<2$. Combining this with Lemma \ref{theorem-spectral}, it suffices to verify that \eqref{estimate-condition} implies \eqref{estimate-condition1}. Specifically, $f(s)=m(s)$, we need to show
\begin{equation}
\sup_{t>0}\|\varphi(\cdot)f(t\cdot)\|_{C^{\sigma}}<\infty,\quad 1<\sigma<2.
\end{equation}
Note that $C^{\sigma}=C^{[\sigma], \sigma-[\sigma]}$. Thus, for $x,y\in\R$ and supp$\varphi\subset[\frac{1}{4},1]$, it is sufficient to prove
\begin{equation}
\sup_{t>0}\left(\|\varphi(\cdot)f(t\cdot)\|_{C^{1}}+\sup_{\substack{x,y\in \textrm{supp}\varphi \\ x\neq y}}\frac{|(\varphi(x)f(tx))'-(\varphi(y)f(ty))'|}{|x-y|^{\sigma-1}}\right)<+\infty.
\end{equation}
We observe that
\begin{equation}\label{estimate-condition2}
|f^{j}(x)|\leq C (1+x)^{-\sigma}, \quad\forall x \geq 0 \text { and } 0 \leq j \leq 2.
\end{equation}
Then, for $1<\sigma<2$, we have
\begin{align}
\sup_{t>0}\|\varphi(\cdot)f(t\cdot)\|_{C^{1}}&\leq\sup_{t>0}|\varphi(x)f(tx)|+\sup_{t>0}|(\varphi(x)f(tx))'|\nonumber\\
&\lesssim\sup_{t>0}(1+t)(1+tx)^{-\sigma}\lesssim 1.
\end{align}

Next, we compute
\begin{align}
\frac{|(\varphi(x)f(tx))'-(\varphi(y)f(ty))'|}{|x-y|^{\sigma-1}}
&=\frac{|\varphi'(x)(f(tx)-f(ty))+(\varphi'(x)-\varphi'(y))f(ty)|}{|x-y|^{\sigma-1}}\nonumber\\
&\quad+\frac{|t\varphi(x)(f'(tx)-f'(ty))+t(\varphi(x)-\varphi(y))f'(ty)|}{|x-y|^{\sigma-1}}\nonumber\\
&\triangleq H_{1}+H_{2}.
\end{align}
\textbf{The proof of $H_{1}$:} For $|x-y|\geq1$, we utilize \eqref{estimate-condition2} in combination with the support of $\varphi(x)$ to establish the boundedness of $H_{1}$. Next, we examine the case where $|x-y|\leq1$.

By applying the mean value theorem, we obtain the following expression
\begin{equation}\label{estimate-condition3}
f(tx)-f(ty)=f'(t\xi_{1})(tx-ty), \quad \xi_{1}\in(\frac{1}{4},1),
\end{equation}
and
\begin{equation}\label{estimate-condition4}
(\varphi'(x)-\varphi'(y))f(ty)=\varphi''(\xi_{2})(x-y)f(ty),\quad \xi_{2}\in(\frac{1}{4},1).
\end{equation}
So, from \eqref{estimate-condition2}, \eqref{estimate-condition3}, and \eqref{estimate-condition4}, as well as the condition $1<\sigma<2$, we derive
\begin{align}
&\sup_{t>0}\frac{|\varphi'(x)(f(tx)-f(ty))+(\varphi'(x)-\varphi'(y))f(ty)|}{|x-y|^{\sigma-1}}\nonumber\\
\lesssim&\sup_{t>0}\frac{t|x-y|^{2-\sigma}}{(1+t\xi_{1})^{\sigma}}+\sup_{t>0}\frac{|x-y|^{2-\sigma}}{(1+t\xi_{2})^{\sigma}}\lesssim 1.
\end{align}
\textbf{The proof of $H_{2}$:} Similarly, for $|tx-ty|\geq1$, the boundedness follows directly from \eqref{estimate-condition2}. We will assume $|tx-ty|\leq1$ for the subsequent analysis.

Using the mean value theorem again, we arrive at
\begin{align}
&\sup_{t>0}\sup_{x\neq y}\frac{|t\varphi(x)(f'(tx)-f'(ty))+t(\varphi(x)-\varphi(y))f'(ty)|}{|x-y|^{\sigma-1}}\nonumber\\
=&\sup_{t>0}\sup_{x\neq y}\frac{|t^{\sigma}\varphi(x)(f'(tx)-f'(ty))+t^{\sigma}(\varphi(x)-\varphi(y))f'(ty)|}{|t(x-y)|^{\sigma-1}}\nonumber\\
\lesssim& t^{\sigma}(1+t)^{-\sigma}|tx-ty|^{2-\sigma}\lesssim 1.
\end{align}
Therefore, we have completed the proof of \eqref{ResultIII}.
\end{proof}

%
%
%
%

\section{The proof of Theorem \ref{operator-family}}\label{proof-family}
In this section, we first construct the kernel of analytic operator family $f_{w,t}(\mathcal{L}_{\mathbf{A}})$, and then we prove the pointwise estimate of the kernel. Finally we show Theorem \ref{operator-family} via the derived pointwise estimate of kernel.

\subsection{The kernel of analytic operator family $
f_{w,t}(\mathcal{L}_{\mathbf{A}})$.} We now draw our attention to setting up the representation of the kernel for analytic operator family $
f_{w,t}(\mathcal{L}_{\mathbf{A}})$ in this subsection.

More precisely, we have the following proposition.
\begin{proposition}[Kernel of the operator $f_{w,t}(\mathcal{L}_{\mathbf{A}})$]\label{prop:reso} Let $x=(r_{1}\cos\theta_{1}, r_{1}\sin\theta_{1})$ and $y=(r_{2}\cos\theta_{2}, r_{2}\sin\theta_{2})$, and we denote $K_{w,t}(r_{1}, \theta_{1}, r_{2}, \theta_{2})$ by the kernel of the operator $f_{w,t}(\mathcal{L}_{\mathbf{A}})$ in \eqref{def:ana-oper}.
Set \begin{align}\label{beta1beta2}
\beta_{1}=\arccos\frac{r_{1}^{2}+r_{2}^{2}-t^{2}}{2r_{1}r_{2}} \quad\text{and}\quad
\beta_{2}=\cosh^{-1}\frac{t^{2}-r_{1}^{2}-r_{2}^{2}}{2r_{1}r_{2}}.
\end{align}
For $t\geq 0$,
then
\begin{equation}
K(t, x, y)=K_{w,t}(r_{1}, \theta_{1}, r_{2},\theta_{2})=G(t,r_{1}, \theta_{1}, r_{2}, \theta_{2})+D(t, r_{1}, \theta_{1}, r_{2},\theta_{2}),
\end{equation}
where
\begin{align}\label{kernelG}
&G(t,r_{1},\theta_{1}, r_{2},\theta_{2})\nonumber\\
=&\frac{t^{2(w-1)}}{\sqrt{2\pi}2^{1-w}\Gamma(1-w)}e^{i\int_{\theta_{2}}^{\theta_{1}}\alpha(\theta')d\theta'}
\left(t^{2}-(r_{1}^{2}+r_{2}^{2}-2r_{1}r_{2}\cos(\theta_{1}-\theta_{2}))\right)^{-w} \nonumber\\
&\times\mathbbm{1}_{(|r_{1}-r_{2}|,r_{1}+r_{2})}(t)[\mathbbm{1}_{[0,\beta_{1}]}(|\theta_{1}-\theta_{2}|)+e^{-2\alpha\pi i} \mathbbm{1}_{[2\pi-\beta_{1},2\pi]}(|\theta_{1}-\theta_{2})],
\end{align}
and
\begin{align}\label{kernelD}
&D(t, r_{1}, \theta_{1}, r_{2}, \theta_{2})\nonumber\\
=&\frac{t^{2(w-1)}e^{i\pi(w-\frac{1}{2})}}{2^{1-w}\pi\sqrt{2\pi}\Gamma(1-w)}\mathbbm{1}_{(r_{1}+r_{2},\infty)}(t)e^{-i(\alpha(\theta_{1}-\theta_{2})-\int_{\theta_{2}}^{\theta_{1}}\alpha(\theta')d\theta')}\nonumber\\
&\times\int_{\beta_{2}}^{+\infty}(t^{2}-r_{1}^{2}-r_{2}^{2}-2r_{1}r_{2}\cosh s)^{-w}B_{\alpha}(s,\theta_{1},\theta_{2})ds,
\end{align}
with
\begin{equation}\label{def:B}
\begin{split}
B_{\alpha}(s,\theta_{1},&\theta_{2}):=\cos(|\alpha|\pi)e^{-|\alpha|s}+\cos(\alpha\pi)\\
&\times\frac{(\cos(\theta_{1}-\theta_{2}+\pi)-e^{-s})\cosh(\alpha s)+i\sin(\theta_{1}-\theta_{2}+\pi)\sinh(\alpha s)}{\cosh s-\cos(\theta_{1}-\theta_{2}+\pi)}.
\end{split}
\end{equation}
When $t\leq 0$, the similar conclusion holds for \eqref{kernelG} and \eqref{kernelD} with replacing $t$ by $-t$.
\end{proposition}

\begin{proof}
With a view to demonstrate the kernel of the operator $f_{w,t}(\mathcal{L}_{\mathbf{A}})$, we start by recalling the kernel of the operator family $f_{w,t}(\mathcal{L}_{\mathbf{A}})$ (see \eqref{FKmuk} and \eqref{Kmuk}). We study the kernel
$$
K_{w,t}(r_{1}, \theta_{1}, r_{2}, \theta_{2})=\sum_{k\in \mathbb{Z}}\varphi_{k}(\theta_{1})\overline{\varphi_{k}(\theta_{2})}K_{\nu_{k}}(t, r_{1}, r_{2}),
$$
where
\begin{align}
K_{\nu_{k}}(t,r_{1},r_{2})&=\int_{0}^{\infty}f_{w,t}(\rho^{2})J_{\nu_{k}}(r_{1}\rho)J_{\nu_{k}}(r_{2}\rho)\rho d\rho\nonumber\\
&=\big(\frac{\pi}{2}\big)^{\frac{1}{2}}\int_{0}^{\infty}(t\rho)^{w-1}J_{1-w}(t\rho)J_{\nu_{k}}(r_{1}\rho)J_{\nu_{k}}(r_{2}\rho)\rho d\rho.
\end{align}
For convenience, in view of \eqref{formula:Ma}, we define
$$\Omega_{I}\triangleq \{(t, r_1, r_2): t<|r_{1}-r_{2}|\},$$
$$\Omega_{II}\triangleq \{(t, r_1, r_2): |r_{1}-r_{2}|<t<r_{1}+r_{2}\},$$
and
$$\Omega_{III}\triangleq \{(t, r_1, r_2): t>r_{1}+r_{2}\}.$$
Then, by  \eqref{formula:Ma}, \eqref{gamma} and \eqref{gammaQ}, we get
\begin{align}\label{kernel-region}
&K_{\nu_{k}}(t, r_{1}, r_{2})\nonumber\\
=& \left(\frac{\pi}{2}\right)^{\frac{1}{2}} t^{w-1} \times
\begin{cases}
0, & \text{if } (t,r_{1},r_{2})\in\Omega_{I}, \\
\frac{(r_{1} r_{2})^{-w} \sin^{\frac{1}{2} - w} \beta_{1}}{(2 \pi)^{\frac{1}{2}} t^{\mathrm{l} - w}} P_{\nu_{k} - \frac{1}{2}}^{w - \frac{1}{2}} (\cos \beta_{1}), & \text{if } (t,r_{1},r_{2})\in\Omega_{II}, \\
\frac{(r_{1} r_{2})^{-w} \sinh^{\frac{1}{2} - w} \beta_{2}}{\left(\frac{\pi^{3}}{2}\right)^{\frac{1}{2}} t^{1 - w}} \cos(\pi \nu_{k}) Q^{w - \frac{1}{2}}_{\nu_{k} - \frac{1}{2}} (\cosh \beta_{2}), & \text{if } (t,r_{1},r_{2})\in\Omega_{III},
\end{cases}\nonumber\\
=& \frac{t^{2(w-1)} (r_{1} r_{2})^{-w}}{\sqrt{2 \pi} \Gamma(1 - w)} \times
\begin{cases}
0, & \text{if } (t,r_{1},r_{2})\in\Omega_{I}, \\
\int_{0}^{\beta_{1}} (\cos s - \cos \beta_{1})^{-w} \cos(\nu_{k} s) \, ds, & \text{if } (t,r_{1},r_{2})\in\Omega_{II}, \\
e^{i \pi (w - \frac{1}{2})} \int_{\beta_{2}}^{+\infty} (\cosh s - \cosh \beta_{2})^{-w} e^{-\nu_{k} s} \cos(\pi \nu_{k}) \, ds, & \text{if } (t,r_{1},r_{2})\in\Omega_{III}.
\end{cases}
\end{align}
Eventually, we can obtain
\begin{align}
K_{\nu_{k}}(t,r_{1},r_{2})
=K_{\nu_{k}}^{\Omega_{I}}(t,r_{1},r_{2})+K_{\nu_{k}}^{\Omega_{II}}(t,r_{1},r_{2})+K_{\nu_{k}}^{\Omega_{III}}(t,r_{1},r_{2}),
\end{align}
where $K_{\nu_k}^{\Omega_{I}}, K_{\nu_{k}}^{\Omega_{II}}$ and $K_{\nu_{k}}^{\Omega_{III}}$ are defined in the region $\Omega_{I}, \Omega_{II}$ and $\Omega_{III}$ respectively. Combining \eqref{kernel-region} with $K_{\nu_{k}}^{\Omega_{I}}\left(t, r_1, r_2\right)=0$, it suffices to calculate the kernel of $K_{\nu_{k}}^{\Omega_{II}}\left(t, r_1, r_2\right)$ and $K_{\nu_{k}}^{\Omega_{III}}\left(t, r_1, r_2\right)$.

Now, we proceed with the kernel of $K_{\nu_{k}}^{\Omega_{II}}\left(t, r_1, r_2\right)$. From the eigenfunction in \eqref{eig-LA}, we need to consider
\begin{align}\label{eig-kernel}
& \sum_{k \in \mathbb{Z}} \frac{1}{2 \pi} e^{-i\left(\left(\theta_1-\theta_2\right)(k+\alpha)-\int_{\theta_2}^{\theta_1} A\left(\theta^{\prime}\right) d \theta^{\prime}\right)}K_{\nu_k}^{\Omega_{II}}\left(t, r_1, r_2\right)\nonumber\\
= & \frac{1}{2 \pi} e^{-i \alpha\left(\theta_1-\theta_2\right)+i \int_{\theta_2}^{\theta_1} A\left(\theta^{\prime}\right) d \theta^{\prime}} \sum_{k \in \mathbb{Z}} e^{-i k\left(\theta_1-\theta_2\right)}K_{\nu_k}^{\Omega_{II}}(t,r_1, r_2).
\end{align}
Firstly, we recall $\nu_k=|k+\alpha|$ with $k \in \mathbb{Z}$ and apply the Poisson summation formula \cite[Theorem 0.1.16]{sogge}, this yields
\begin{align}\label{sum-cos}
\sum_{k \in \mathbb{Z}}\frac{1}{2\pi} \cos \left(s \nu_k\right) e^{-i k\left(\theta_1-\theta_2\right)}
=\frac{1}{2}\sum_{k\in \Z}(e^{i\alpha s}\delta(\theta_{1}-\theta_{2}-s+2\pi k)+e^{-i\alpha s}\delta(\theta_{1}-\theta_{2}+s+2\pi k)).
\end{align}

To facilitate the calculation, we will use $\bar{\theta}$ to represent $\theta_{1}-\theta_{2}$ in the subsequent discussion. Thus, if $\left|r_1-r_2\right|<t<r_1+r_2$, from \eqref{beta1beta2}, \eqref{kernel-region} and \eqref{sum-cos}, we can derive the following expression
$$
\begin{aligned}
& \sum_{k \in \mathbb{Z}} e^{-i k\bar{\theta}} K_{\nu_{k}}^{\Omega_{II}}\left(t, r_1, r_2\right) \\
= & \frac{t^{2(w-1)}\left(r_1 r_2\right)^{-w}}{\sqrt{2 \pi} \Gamma(1-w)} \sum_{k \in \mathbb{Z}} e^{-i k\bar{\theta}} \int_0^{\beta_1}\left(\cos s-\cos \beta_1\right)^{-w} \cos \left(\nu_k s\right) d s \\
= & \frac{\sqrt{2 \pi}t^{2(w-1)}\left(r_1 r_2\right)^{-w}}{2 \Gamma(1-w)} \sum_{k \in \mathbb{Z}} \int_0^{\beta_1} \frac{e^{-i s \alpha} \delta\left(\bar{\theta}+2 k \pi+s\right)+e^{i s \alpha} \delta\left(\bar{\theta}+2 k \pi-s\right)}{\left(\cos s-\cos \beta_1\right)^w} d s \\
= & \frac{ \sqrt{2 \pi}t^{2(w-1)}}{2^{1-w} \Gamma(1-w)}\times\sum_{\left\{k \in \mathbb{Z}: 0 \leq\left|\bar{\theta}+2 k \pi\right| \leq \beta_1\right\}}\left[t^2-r_1^2-r_2^2+2 r_1 r_2 \cos \left(\bar{\theta}+2 k \pi\right)\right]^{-w} e^{i\left(\bar{\theta}+2 k \pi\right) \alpha} \\
= & \frac{ \sqrt{2 \pi}t^{2(w-1)}}{2^{1-w} \Gamma(1-w)} \times \begin{cases}0, & \text { if } \beta_1<\left|\bar{\theta}\right|<2 \pi-\beta_1, \\
{\left[t^2-r_1^2-r_2^2+2 r_1 r_2 \cos (\bar{\theta})\right]^{-w} e^{i\bar{\theta} \alpha}}, & \text { if }\left|\bar{\theta}\right|<\beta_1, \\
{\left[t^2-r_1^2-r_2^2+2 r_1 r_2 \cos \left(\bar{\theta}\right)\right]^{-w} e^{i\left(\bar{\theta}-2 \pi\right) \alpha}}, & \text { if } 2\pi-\beta_1<\left|\bar{\theta}\right|<2 \pi .\end{cases}
\end{aligned}
$$
Hence, by \eqref{eig-kernel}, we obtain the term in $G(t,r_{1},r_{2},\theta_{1},\theta_{2})$ as follows
\begin{align}\label{Gkernel}
& \frac{t^{2(w-1)}}{\sqrt{2 \pi} 2^{1-w} \Gamma(1-w)}\left(t^2-\left(r_1^2+r_2^2-2 r_1 r_2 \cos \left(\bar{\theta}\right)\right)\right)^{-w} e^{i \int_{\theta_2}^{\theta_1} \alpha\left(\theta^{\prime}\right) d \theta^{\prime}} \nonumber\\
& \quad\times \mathbbm{1}_{\left(\left|r_1-r_2\right|, r_1+r_2\right)}(t)\left[\mathbbm{1}_{\left[0, \beta_1\right]}(|\bar{\theta}|)+e^{-2 \alpha \pi i} \mathbbm{1}_{\left[2 \pi-\beta_1, 2 \pi\right]}(|\bar{\theta}|)\right] .
\end{align}

Next, we take account of the kernel function
\begin{align*}
K_{\nu_{k}}^{\Omega_{III}}\left(t, r_1, r_2\right)=\frac{t^{2(w-1)}e^{i\pi(w-\frac{1}{2})}(r_{1}r_{2})^{-w}}{\sqrt{2 \pi} \Gamma(1-w)}\int_{\beta_{2}}^{+\infty}(\cosh s-\cosh\beta_{2})^{-w}e^{-\nu_{k}s}\cos(\pi\nu_{k})ds
\end{align*}
in region $t>r_{1}+r_{2}$.
Thus, it remains to evaluate the integral
\begin{equation}
\int_{\beta_{2}}^{+\infty}(\cosh s-\cosh\beta_{2})^{-w}e^{-\nu_{k}s}\cos(\pi\nu_{k})ds.
\end{equation}

Due to the unitary equivalence, it is enough to think about the case $\alpha \in(0,1)$, then
$$
\nu_{k}=|k+\alpha|= \begin{cases}k+\alpha, & k \geq 1, \\ |\alpha|, & k=0, \\ -(k+\alpha), & k \leq-1 .\end{cases}
$$
Therefore, similar to the process of \eqref{eig-kernel}, we furthermore compute
$$
\begin{aligned}
& \sum_{k \in \mathbb{Z}} \cos (\pi|k+\alpha|) e^{-s|k+\alpha|} e^{-i k\bar{\theta}} \\
= & \cos (\alpha \pi) \sum_{k \geq 1} \frac{e^{i k \pi}+e^{-i k \pi}}{2} e^{-s(k+\alpha)} e^{-i k\bar{\theta}}+\cos (|\alpha| \pi) e^{-|\alpha| s} \\
& +\cos (\alpha \pi) \sum_{k \leq-1} \frac{e^{i k \pi}+e^{-i k \pi}}{2} e^{s(k+\alpha)} e^{-i k\bar{\theta}} \\
= & \cos (|\alpha| \pi) e^{-|\alpha| s}+\frac{\cos (\alpha \pi)}{2}\Big(e^{-s \alpha} \sum_{k \geq 1} e^{-k s}(e^{-i k(\bar{\theta}+\pi)}+e^{-i k(\bar{\theta}-\pi)})\Big. \\
&\Big.+e^{s \alpha} \sum_{k \geq 1} e^{-k s}(e^{i k(\bar{\theta}+\pi)}+e^{i k(\bar{\theta}-\pi)})\Big) .
\end{aligned}
$$

Observe that
$$
\sum_{k=1}^{\infty} e^{i k z}=\frac{e^{i z}}{1-e^{i z}}, \quad \operatorname{Im} z>0,
$$
we finally achieve
$$
\begin{aligned}
& \sum_{k \in \mathbb{Z}} \cos (\pi|k+\alpha|) e^{-s|k+\alpha|} e^{-i k\bar{\theta}} \\
=&\cos (|\alpha| \pi) e^{-|\alpha| s}+\frac{\cos (\alpha \pi)}{2}\left(\frac{e^{-(1+\alpha) s-i\left(\bar{\theta}+\pi\right)}}{1-e^{-s-i\left(\bar{\theta}+\pi\right)}}+\frac{e^{-(1+\alpha) s-i\left(\bar{\theta}-\pi\right)}}{1-e^{-s-i\left(\bar{\theta}-\pi\right)}}\right. \\
&\left.\quad+\frac{e^{-(1-\alpha) s+i\left(\bar{\theta}+\pi\right)}}{1-e^{-s+i\left(\bar{\theta}+\pi\right)}}+\frac{e^{-(1-\alpha) s+i\left(\bar{\theta}-\pi\right)}}{1-e^{-s+i\left(\bar{\theta}-\pi\right)}}\right) \\
=&\cos (|\alpha| \pi) e^{-|\alpha| s}+\cos (\alpha \pi) \frac{\left(\cos \left(\bar{\theta}+\pi\right)-e^{-s}\right) \cosh (\alpha s)+i \sin \left(\bar{\theta}+\pi\right) \sinh (\alpha s)}{\cosh (s)-\cos \left(\bar{\theta}+\pi\right)}, \\
&
\end{aligned}
$$
which gives the term $D(t,r_{1},r_{2},\theta_{1},\theta_{2})$.

\end{proof}
\subsection{Pointwise estimate for $G(r_{1}, r_{2}, \theta_{1}, \theta_{2})$ and $D(r_{1}, r_{2}, \theta_{1}, \theta_{2})$}
In this subsection, we mainly explore to establish the pointwise estimate for the terms $G(r_{1}, r_{2}, \theta_{1}, \theta_{2})$ and $D(r_{1}, r_{2}, \theta_{1}, \theta_{2})$, which is crucial for the establishment of Theorem \ref{operator-family}.

Our main result of this subsection is the following Proposition.
\begin{proposition}[Pointwise estimate for the kernel of $f_{w,t}(\mathcal{L}_{\mathbf{A},a})$]\label{pointwise}
Assume that $\mathrm{Re} w\in(0,1)$. Let $G(r_{1}, r_{2}, \theta_{1}, \theta_{2})$ be as in \eqref{kernelG} and $D(r_{1}, r_{2}, \theta_{1}, \theta_{2})$ be as in \eqref{kernelD}. Then, the following estimates hold
\begin{align}\label{estimate:G1}
|G(r_{1}, r_{2}, \theta_{1}, \theta_{2})|\lesssim \frac{t^{2(\mathrm{Re} w-1)}}{(t^2-|x-y|^2)^{\mathrm{Re} w}}, \quad (r_{1}+r_{2})^2>t^2>|x-y|^2,
\end{align}

\begin{align}\label{estimate:D}
|D(r_{1}, r_{2}, \theta_{1}, \theta_{2})|\lesssim \frac{t^{2(\mathrm{Re} w-1)}}{(t^2-(r_{1}+r_{2})^2)^{\mathrm{Re} w}}, \quad t^2>(r_{1}+r_{2})^2.
\end{align}
\end{proposition}
\begin{proof}
First of all, in the polar coordinates $x=r_1(\cos\theta_1,\sin\theta_1)$ and $y=r_{2}(\cos\theta_{2},\sin\theta_2)$, the squared distance between the two points can be calculated as
\begin{equation}
|x-y|^2=r_{1}^2+r_{2}^2-2r_{1}r_{2}\cos(\theta_{1}-\theta_{2}).
\end{equation}
From this, we can derive the inequality $t^2>(r_{1}-r_{2})^2.$
Thus, we can conclude that the two sets
\begin{equation*}
\{|\theta_{1}-\theta_{2}|:t^2>r_{1}^2+r_{2}^2-2r_{1}r_{2}\cos(\theta_{1}-\theta_{2})\}\Longleftrightarrow \{|\theta_{1}-\theta_{2}|\in [0,\beta_{1}]\cup[2\pi-\beta_{1},2\pi]\}
\end{equation*}
given the conditions $|r_{1}-r_{2}|<t<r_{1}+r_{2}$ and $|\theta_{1}-\theta_{2}|\leq 2\pi$. Indeed, this relationship can be verified based on the condition that $\beta_{1}\in[0,\pi]$, as specified in \eqref{beta1beta2}.

Next, from the expression of the term $G(t, r_1, \theta_1, r_2, \theta_2)$ in \eqref{kernelG}, we can deduce
$$
\begin{aligned}
&\left|G(t, r_1, \theta_1, r_2, \theta_2)\right| \\
\leq&\Big|\frac{t^{2(w-1)}}{\sqrt{2 \pi} 2^{1-w} \Gamma(1-w)}\Big|\left(t^2-\left(r_1^2+r_2^2-2 r_1 r_2 \cos (\bar{\theta})\right)\right)^{-\mathrm{Re} w} \\
&\quad \times\mathbbm{1}_{\left(\left|r_1-r_2\right|, r_1+r_2\right)}(t)[\mathbbm{1}_{\left[0, \beta_1\right]}(|\bar{\theta}|)+\mathbbm{1}_{\left[2 \pi-\beta_1, 2 \pi\right]}(|\bar{\theta}|)] \\
\lesssim& t^{2(\mathrm{Re} w-1)}\left(t^2-|x-y|^2\right)^{-\mathrm{Re} w},
\end{aligned}
$$
hence we prove \eqref{estimate:G1}.

Now, we verify \eqref{estimate:D}. To this end, from the definition of $B_{\alpha}(s,\theta_{1},\theta_{2})$ in \eqref{def:B}, we are required to estimate the following three terms:
\begin{align}\label{estimate:D1}
\cos (|\alpha| \pi) \int_{\beta_2}^{+\infty} \left(t^2-r_1^2-r_2^2-2r_1r_2\cosh \tau\right)^{-w} e^{-|\alpha| \tau} d\tau,
\end{align}
\begin{align}\label{estimate:D2}
\cos (\alpha \pi) \int_{\beta_2}^{+\infty} \left(t^2-r_1^2-r_2^2-2r_1r_2\cosh \tau\right)^{-w} \frac{\left(\cos \left(\bar{\theta} + \pi\right)-e^{-\tau}\right) \cosh (\alpha \tau)}{\cosh (\tau) - \cos \left(\bar{\theta} + \pi\right)} d\tau,
\end{align}
\begin{align}\label{estimate:D3}
\cos (\alpha \pi) \int_{\beta_2}^{+\infty} \left(t^2-r_1^2-r_2^2-2r_1r_2\cosh \tau\right)^{-w} \frac{\sin \left(\bar{\theta} + \pi\right) \sinh (\alpha \tau)}{\cosh (\tau) - \cos \left(\bar{\theta} + \pi\right)} d\tau.
\end{align}

\textbf{Contribution of \eqref{estimate:D1}:} We first consider \eqref{estimate:D1}. The following simple calculation indicates
\begin{align}\label{estimate-calculus}
\frac{1}{(\cosh (\beta_2)-\cosh \tau)^{w}}&=\frac{1}{(2 \sinh \frac{\beta_2+\tau}{2} \sinh \frac{\beta_2-\tau}{2})^w}\nonumber\\
&=\frac{\sinh ^{2 w}(\frac{\beta_2}{2})}{\sinh ^{2 w}(\frac{\beta_2}{2})(2 \sinh \frac{\beta_2+\tau}{2} \sinh \frac{\beta_2-\tau}{2})^w}\nonumber\\
&=\frac{2^w(2r_{1}r_{2})^w\sinh ^{2 w}(\frac{\beta_2}{2})}{(t^2-(r_1+r_{2})^2)^{w}(2 \sinh \frac{\beta_2+\tau}{2} \sinh \frac{\beta_2-\tau}{2})^w},
\end{align}
where we have used
$$
\cosh \alpha-\cosh \beta=2 \sinh \frac{\alpha+\beta}{2} \sinh \frac{\alpha-\beta}{2},
$$
and

\begin{align}\label{estimate:sinbeta2}
\sinh \left(\frac{\beta_2}{2}\right) & =\sqrt{\cosh ^2\left(\frac{\beta_2}{2}\right)-1}=\sqrt{\frac{\cosh \left(\beta_2\right)-1}{2}} \nonumber\\
& =\frac{1}{\sqrt{2}} \sqrt{\frac{t^2-r_1^2-r_2^2}{2 r_1 r_2}-1}=\frac{1}{\sqrt{2}} \sqrt{\frac{t^2-\left(r_1+r_2\right)^2}{2 r_1 r_2}} .
\end{align}
Therefore, we can obtain

\begin{align}\label{estimate:e}
& \int_{\beta_2}^{+\infty}\left(t^2-r_1^2-r_2^2-2 r_1 r_2 \cosh \tau\right)^{-w} e^{-|\alpha| \tau} d \tau \nonumber\\
= & \frac{1}{\left(2 r_1 r_2\right)^w} \int_{\beta_2}^{+\infty}\big(\frac{t^2-r_1^2-r_2^2}{2 r_1 r_2}-\cosh \tau\big)^{-w} e^{-|\alpha| \tau} d \tau \nonumber\\
= & \frac{1}{\left(2 r_1 r_2\right)^w} \int_{\beta_2}^{+\infty}\left(\cosh (\beta_2)-\cosh \tau\right)^{-w} e^{-|\alpha| \tau} d \tau \nonumber\\
= & \frac{2^w}{\big(t^2-\left(r_1+r_2\right)^2\big)^w} \int_{\beta_2}^{+\infty} \frac{\sinh { }^{2 w}\left(\frac{\beta_2}{2}\right)}{e^{|\alpha| \tau}\big(2 \sinh \frac{\beta_2+\tau}{2} \sinh \frac{\beta_2-\tau}{2}\big)^w} d \tau.
\end{align}
Now we claim that
$$
\Big|\int_{\beta_2}^{+\infty} \frac{\sinh ^{2 w}\left(\frac{\beta_2}{2}\right)}{e^{|\alpha| \tau}\big(2 \sinh \frac{\beta_2+\tau}{2} \sinh \frac{\beta_2-\tau}{2}\big)^w} d \tau\Big| \lesssim 1 .
$$

\textbf{Case 1: $\beta_2 \geq 1$.} By computation, we gain
\begin{align}\label{estimate-beta2}
& \left|\int_{\beta_2}^{+\infty} \frac{\sinh ^{2 w}\left(\frac{\beta_2}{2}\right)}{e^{|\alpha| \tau}\big(2 \sinh \frac{\beta_2+\tau}{2} \sinh \frac{\beta_2-\tau}{2}\big)^ w} d \tau\right| \nonumber\\
\lesssim &\int_{\beta_2}^{+\infty} \frac{e^{\beta_2 \mathrm{Re} w}}{e^{|\alpha| \tau}\big(e^{\frac{\beta_2+\tau}{2}}-e^{-\frac{\beta_2+\tau}{2}}\big)^{\mathrm{Re}w}\big(e^{\frac{\tau-\beta_2}{2}}-e^{-\frac{\tau-\beta_2}{2}}\big)^{\mathrm{Re} w}} d \tau \nonumber\\
\lesssim& \int_{\beta_2}^{+\infty} \frac{e^{\beta_2 \mathrm{Re} w}}{e^{|\alpha| \tau} e^{\frac{\beta_2+\tau}{2} \mathrm{Re} w}\big(e^{\frac{\tau-\beta_2}{2}}-e^{-\frac{\tau-\beta_2}{2}}\big)^{\mathrm{Re} w}} d \tau \nonumber\\
=&\int_{\beta_2}^{+\infty} \frac{1}{e^{|\alpha| \tau}(e^{\tau-\beta_2}-1)^{\mathrm{Re} w}} d \tau \nonumber\\
=&\int_{\beta_2}^{\beta_2+\frac{1}{2}} \frac{1}{e^{|\alpha| \tau}(e^{\tau-\beta_2}-1)^{\mathrm{Re} w}} d \tau+\int_{\beta_2+\frac{1}{2}}^{\infty} \frac{1}{e^{|\alpha| \tau}(e^{\tau-\beta_2}-1)^{\mathrm{Re} w}} d \tau \nonumber\\
\lesssim& \int_{\beta_2}^{\beta_2+\frac{1}{2}} \frac{1}{e^{|\alpha| \tau}\left(\tau-\beta_2\right)^{\mathrm{Re} w}} d \tau+\int_{\beta_2+\frac{1}{2}}^{+\infty} e^{-|\alpha| \tau} d \tau \lesssim 1 . \nonumber\\
&
\end{align}

\textbf{Case 2: $\beta_2 \leq 1$.} In this case, we can acquire
\begin{align}\label{estimate:beta2}
& \left|\int_{\beta_2}^{+\infty} \frac{\sinh ^{2 w}\left(\frac{\beta_2}{2}\right)}{e^{|\alpha| \tau}\big(2 \sinh \frac{\beta_2+\tau}{2} \sinh \frac{\beta_2-\tau}{2}\big)^w} d \tau\right| \nonumber\\
\lesssim & \int_{\beta_2}^{\beta_2+1} \frac{e^{-\alpha \tau}\beta_2^{2 \mathrm{Re} w}}{\left(\beta_2+\tau\right)^{\mathrm{Re} w}(\tau-\beta_2)^{\mathrm{Re} w}} d \tau+\int_{\beta_2+1}^{+\infty}\frac{e^{\beta_2 \mathrm{Re} w}}{e^{|\alpha| \tau} e^{\frac{\beta_2+\tau}{2} \mathrm{Re} w}\big(e^{\frac{\tau-\beta_2}{2}}-e^{-\frac{\tau-\beta_2}{2}}\big)^{\mathrm{Re} w}} d \tau  \nonumber\\
\lesssim&\int_{\beta_2}^{\beta_2+1} \frac{\beta_2^{\mathrm{Re} w}}{(\tau-\beta_2)^{\mathrm{Re} w}} d \tau + \int_{\beta_2+1}^{+\infty}\frac{1}{e^{|\alpha| \tau}(e^{\tau-\beta_2}-1)^{\mathrm{Re} w}} d \tau \nonumber\\
\lesssim&\int_{\beta_2}^{\beta_2+1} \frac{1}{(\tau-\beta_2)^{\mathrm{Re} w}} d \tau+\int_{\beta_2+1}^{+\infty}e^{-\alpha\tau}d\tau \lesssim 1,
\end{align}
where we have used $\mathrm{Re} w<1$. This together with \eqref{estimate-beta2} and \eqref{estimate:beta2} yields that
$$
\int_{\beta_2}^{+\infty}\left(t^2-r_1^2-r_2^2-2 r_1 r_2 \cosh \tau\right)^{-w} e^{-|\alpha| \tau} d \tau \lesssim \frac{1}{(t^2-\left(r_1+r_2\right)^2)^{\operatorname{Re} w}} .
$$

\textbf{Contribution of \eqref{estimate:D2}:} Note that
$$
\cosh (\tau)-\cos \left(\bar{\theta}+\pi\right)=\sinh ^2\left(\frac{\tau}{2}\right)+\sin ^2\left(\frac{\bar{\theta}+\pi}{2}\right),
$$
and by \eqref{estimate-calculus}, we are going to evaluate
\begin{align}\label{D2sin}
&\int_{\beta_2}^{+\infty}\left(t^2-r_1^2-r_2^2-2 r_1 r_2 \cosh \tau\right)^{-w} \frac{\cosh (\alpha \tau)\left(\cos \left(\bar{\theta}+\pi\right)-e^{-\tau}\right)}{\sinh ^2\left(\frac{\tau}{2}\right)+\sin ^2\left(\frac{\bar{\theta}+\pi}{2}\right)} d \tau \nonumber\\
= & \frac{2^w}{\left(t^2-\left(r_1+r_2\right)^2\right)^w} \int_{\beta_2}^{+\infty} \frac{\sinh ^{2 w}\left(\frac{\beta_2}{2}\right)}{\left(2 \sinh \frac{\beta_2+\tau}{2} \sinh \frac{\beta_2-\tau}{2}\right)^w} \frac{\cosh (\alpha \tau)\left(\cos \left(\bar{\theta}+\pi\right)-e^{-\tau}\right)}{\sinh ^2\left(\frac{\tau}{2}\right)+\sin ^2\left(\frac{\bar{\theta}+\pi}{2}\right)} d \tau \nonumber\\
= & \frac{1}{\left(t^2-\left(r_1+r_2\right)^2\right)^w} \int_{\beta_2}^{+\infty} \frac{\sinh ^{2 w}\left(\frac{\beta_2}{2}\right)}{\left(2 \sinh \frac{\beta_2+\tau}{2} \sinh \frac{\beta_2-\tau}{2}\right)^w} \frac{\cosh (\alpha \tau)\left(\cos \left(\bar{\theta}+\pi\right)-1+1-e^{-\tau}\right)}{\sinh ^2\left(\frac{\tau}{2}\right)+\sin ^2\left(\frac{\bar{\theta}+\pi}{2}\right)} d \tau \nonumber\\
= & \frac{1}{\left(t^2-\left(r_1+r_2\right)^2\right)^w} \int_{\beta_2}^{+\infty} \frac{\sinh ^{2 w}\left(\frac{\beta_2}{2}\right)}{\left(2 \sinh \frac{\beta_2+\tau}{2} \sinh \frac{\beta_2-\tau}{2}\right)^w} \frac{\cosh (\alpha \tau)\left(1-e^{-\tau}-2 \sin ^2\left(\frac{\bar{\theta}+\pi}{2}\right)\right)}{\sinh ^2\left(\frac{\tau}{2}\right)+\sin ^2\left(\frac{\bar{\theta}+\pi}{2}\right)} d \tau .
\end{align}
To estimate the above equality, similarly as above, we divide into two cases.

\textbf{Case 1: $\beta_2 \leq 1$.} On the one hand, from the fact $\alpha<1$, we have

\begin{align}\label{D2sinh}
&\left|\int_{\beta_2}^{+\infty} \frac{\sinh ^{2 w}\left(\frac{\beta_2}{2}\right)}{\left(2 \sinh \frac{\beta_2+\tau}{2} \sinh \frac{\beta_2-\tau}{2}\right)^w} \frac{\cosh (\alpha \tau)\left(1-e^{-\tau}\right)}{\sinh ^2\left(\frac{\tau}{2}\right)+\sin ^2\left(\frac{\bar{\theta}+\pi}{2}\right)} d \tau\right| \nonumber\\
\lesssim& \int_{\beta_2}^{2\beta_2} \frac{\beta_2^{\mathrm{Re} w}}{\left(\tau-\beta_2\right)^{\mathrm{Re} w}} \frac{\tau}{\tau^2+\left(\frac{\bar{\theta}+\pi}{2}\right)^2} d \tau +\int_{2\beta_2}^{+\infty} \frac{e^{\alpha \tau}}{e^{\tau}+\left(\frac{\bar{\theta}+\pi}{2}\right)^2} d \tau \nonumber\\
\lesssim& \frac{1}{\beta_2}\int_{\beta_2}^{2\beta_2} \frac{\beta_2^{\mathrm{Re} w}}{\left(\tau-\beta_2\right)^{\mathrm{Re} w}} d \tau+\int_{2\beta_{2}}^{2\beta_{2}+1}\frac{\beta_2^{\mathrm{Re} w}}{\left(\tau-\beta_2\right)^{\mathrm{Re} w}}\frac{1}{\tau}d \tau +\int_{2\beta_2+1}^{+\infty} e^{-(1-\alpha)\tau}d\tau \lesssim 1,
\end{align}
provided that
\begin{align}
\int_{2\beta_{2}}^{2\beta_{2}+1}\frac{\beta_2^{\mathrm{Re} w}}{\left(\tau-\beta_2\right)^{\mathrm{Re} w}}\frac{1}{\tau}d \tau
&=\int_{2}^{2+\frac{1}{\beta_{2}}}\frac{1}{s(s-1)^{\mathrm{Re} w}}ds\nonumber\\
&\leq \int_{2}^{+\infty}\frac{1}{s(s-1)^{\mathrm{Re} w}}ds\lesssim 1,
\end{align}
in which we use the variable substitution $s=\frac{\tau}{\beta_{2}}$.

On the other hand, similarly, we acquire
\begin{align}\label{D2sinhA}
&\left|\int_{\beta_2}^{+\infty}  \frac{\sinh ^{2 w}\left(\frac{\beta_2}{2}\right)}{\left(2 \sinh \frac{\beta_2+\tau}{2} \sinh \frac{\beta_2-\tau}{2}\right)^w} \frac{\cosh (\alpha \tau) \sin ^2\left(\frac{\bar{\theta}+\pi}{2}\right)}{\sinh ^2\left(\frac{\tau}{2}\right)+\sin ^2\left(\frac{\bar{\theta}+\pi}{2}\right)} d \tau\right| \nonumber\\
\lesssim &\int_{\beta_2}^{\beta_2+1} \frac{\beta_2^{\mathrm{Re} w}}{\left(\tau-\beta_2\right)^{\mathrm{Re} w}} \frac{\left(\frac{\bar{\theta}+\pi}{2}\right)^2}{\tau^2+\left(\frac{\bar{\theta}+\pi}{2}\right)^2} d \tau+\int_{\beta_2+1}^{+\infty}e^{-(1-\alpha)\tau}d\tau \nonumber\\
\lesssim &\int_{\beta_2}^{\beta_2+1} \frac{\beta_2^{\mathrm{Re} w}}{\left(\tau-\beta_2\right)^{\mathrm{Re} w}} d \tau+\int_{\beta_2+1}^{+\infty}e^{-(1-\alpha)\tau}d\tau \lesssim 1.
\end{align}
This together with \eqref{D2sin}-\eqref{D2sinhA} implies
$$
\eqref{estimate:D2} \lesssim \frac{1}{(t^2-\left(r_1+r_2\right)^2)^{\mathrm{Re} w}}
$$
in the case that $\beta_2 \leq 1$.

\textbf{Case 2: $\beta_2 \geq 1$.} By the above argument as \eqref{estimate:D1}, it suffices to take account of the integral in \eqref{estimate:D2}. Then we have
\begin{align}\label{D3sinh}
&\left|\int_{\beta_2}^{+\infty} \frac{\sinh ^{2 w}\left(\frac{\beta_2}{2}\right)}{\left(2 \sinh \frac{\beta_2+\tau}{2} \sinh \frac{\beta_2-\tau}{2}\right)^w} \frac{\cosh (\alpha \tau)\left|\left(e^{-\tau}-1+2 \sin ^2\left(\frac{\bar{\theta}+\pi}{2}\right)\right)\right|}{\sinh ^2\left(\frac{\tau}{2}\right)+\sin ^2\left(\frac{\bar{\theta}+\pi}{2}\right)} d \tau\right| \nonumber\\
\lesssim& \int_{\beta_2}^{+\infty} \frac{e^{\beta_2 \mathrm{Re} w}}{\left[e^{\frac{\beta_2+\tau}{2}}\left(e^{\frac{\tau-\beta_2}{2}}-e^{-\frac{\tau-\beta_2}{2}}\right)\right]^{\mathrm{Re} w}} \frac{\cosh (\alpha \tau)}{\sinh ^2\left(\frac{\tau}{2}\right)+\sin ^2\left(\frac{\bar{\theta}+\pi}{2}\right)} d \tau \nonumber\\
\lesssim& \int_{\beta_2}^{+\infty} \frac{1}{\left(e^{\tau-\beta_2}-1\right)^{\mathrm{Re} w}} \frac{\cosh (\alpha \tau)}{\sinh ^2\left(\frac{\tau}{2}\right)+\sin ^2\left(\frac{\bar{\theta}+\pi}{2}\right)} d \tau \nonumber\\
\lesssim&\int_{\beta_2+\frac{1}{2}}^{+\infty} \frac{\cosh (\alpha \tau)}{\sinh ^2\left(\frac{\tau}{2}\right)+\sin ^2\left(\frac{\bar{\theta}+\pi}{2}\right)} d \tau+\int_{\beta_2}^{\beta_2+\frac{1}{2}} \frac{1}{\left(\tau-\beta_2\right)^{\mathrm{Re} w}} \frac{\cosh (\alpha \tau)}{\sinh ^2\left(\frac{\tau}{2}\right)+\sin ^2\left(\frac{\bar{\theta}+\pi}{2}\right)} d \tau \nonumber\\
\lesssim&\int_{\beta_2+\frac{1}{2}}^{+\infty}  e^{-(1-\alpha) \tau} d \tau+\int_{\beta_2}^{\beta_2+\frac{1}{2}}  \frac{1}{\left(\tau-\beta_2\right)^{\mathrm{Re} w}} e^{-(1-\alpha) \tau} d \tau\nonumber\\
\lesssim &1.
\end{align}
Therefore, we conclude
\begin{align}
\eqref{estimate:D2} \lesssim \frac{1}{(t^2-\left(r_1+r_2\right)^2)^{\mathrm{Re} w}}
\end{align}
in the case that $\beta_2 \geq 1$.

\textbf{Contribution of \eqref{estimate:D3}.} Similarly as \eqref{estimate:D2}, we write
\begin{align}\label{D3write}
& \int_{\beta_2}^{+\infty}\left(t^2-r_1^2-r_2^2-2 r_1 r_2 \cosh \tau\right)^{-w} \frac{\sinh (\alpha \tau) \sin \left(\bar{\theta}+\pi\right)}{\sinh ^2\left(\frac{\tau}{2}\right)+\sin ^2\left(\frac{\bar{\theta}+\pi}{2}\right)} d \tau \nonumber\\
= & \frac{1}{\big(t^2-\left(r_1+r_2\right)^2\big)^w} \int_{\beta_2}^{+\infty} \frac{\sinh ^{2 w}\left(\frac{\beta_2}{2}\right)}{\left(2 \sinh \frac{\beta_2+\tau}{2} \sinh \frac{\beta_2-\tau}{2}\right)^w} \frac{\sinh (\alpha\tau) \sin \left(\bar{\theta}+\pi\right)}{\sinh ^2\left(\frac{\tau}{2}\right)+\sin ^2\left(\frac{\bar{\theta}+\pi}{2}\right)} d \tau .
\end{align}
Next, we divide into two cases.

\textbf{Case 1: $\beta_2 \leq 1$.} Let $b:=|\sin \left(\frac{\bar{\theta}+\pi}{2}\right)|$, we obtain
\begin{align*}
& \left|\int_{\beta_2}^{+\infty} \frac{\sinh ^{2 w}\left(\frac{\beta_2}{2}\right)}{\left(2 \sinh \frac{\beta_2+\tau}{2} \sinh \frac{\beta_2-\tau}{2}\right)^w} \frac{\sinh (\alpha \tau) \sin \left(\bar{\theta}+\pi\right)}{\sinh ^2\left(\frac{\tau}{2}\right)+\sin ^2\left(\frac{\bar{\theta}+\pi}{2}\right)} d \tau\right| \\
\lesssim &\int_{\beta_2}^{\beta_2+1} \frac{\beta_2^{\mathrm{Re} w}}{\left(\tau-\beta_2\right)^{\mathrm{Re} w}} \frac{(\alpha \tau)b}{\tau^2+b^2} d \tau+\int_{\beta_2+1}^{+\infty}e^{-(1-\alpha)\tau}d\tau  \\
\lesssim&\int_{\beta_2}^{\beta_2+1} \frac{\beta_2^{\mathrm{Re} w}}{\left(\tau-\beta_2\right)^{\mathrm{Re} w}} d \tau+\int_{\beta_2+1}^{+\infty}e^{-(1-\alpha)\tau}d\tau \\
\lesssim &1.
\end{align*}

This implies
$$
\eqref{estimate:D3} \lesssim \frac{1}{(t^2-\left(r_1+r_2\right)^2)^{\mathrm{Re} w}}
$$
in the case that $\beta_2 \leq 1$.

\textbf{Case 2: $\beta_2 \geq 1$.} By the above argument and the similar calculations, we can obtain
$$
\begin{aligned}
& \left|\int_{\beta_2}^{+\infty} \frac{\sinh ^{2 w}\left(\frac{\beta_2}{2}\right)}{\left(2 \sinh \frac{\beta_2+\tau}{2} \sinh \frac{\beta_2-\tau}{2}\right)^w} \frac{\sinh (\alpha \tau) \sin \left(\bar{\theta}+\pi\right)}{\sinh ^2\left(\frac{\tau}{2}\right)+\sin ^2\left(\frac{\bar{\theta}+\pi}{2}\right)} d \tau\right| \\
\lesssim &\int_{\beta_2}^{+\infty} \frac{e^{\beta_2 \mathrm{Re} w}}{\left[e^{\frac{\beta_2+\tau}{2}}\left(e^{\frac{\tau-\beta_2}{2}}-e^{-\frac{\tau-\beta_2}{2}}\right)\right]^{\mathrm{Re} w}} \frac{\sinh (\alpha \tau)}{\sinh ^2\left(\frac{\tau}{2}\right)} d \tau \\
\lesssim& \int_{\beta_2}^{+\infty} \frac{1}{\left(e^{\tau-\beta_2}-1\right)^{\mathrm{Re} w}} e^{-(1-\alpha) \tau} d \tau \\
\lesssim&\int_{\beta_2}^{\beta_2+1}\frac{1}{\left(\tau-\beta_2\right)^{\mathrm{Re} w}}d \tau+\int_{\beta_2+1}^{+\infty} e^{-(1-\alpha) \tau} d \tau \\
\lesssim& 1 .
\end{aligned}
$$

This together with \eqref{D3write} yields
$$
\eqref{estimate:D3} \lesssim \frac{1}{(t^2-\left(r_1+r_2\right)^2)^{\mathrm{Re} w}}
$$
in the case where $\beta_2 \geq 1$.

Consequently, we have completed the proof of Proposition \ref{pointwise}.

\end{proof}
\subsection{The proof of Theorem \ref{operator-family}}
In this subsection, we aim to verify Theorem \ref{operator-family} via making use of the pointwise estimate for the kernel of analytic operator $f_{w,t}(\mathcal{L}_{\mathbf{A}})$. According to Proposition \ref{pointwise}, for $1\leq p\leq+\infty$, we can see
\begin{align}\label{shurtest}
& \left\| f_{w, t}\left(\mathcal{L}_{\mathbf{A}}\right) g \right\|_{L^p\left(\mathbb{R}^2\right)} \nonumber \\
\lesssim& t^{2(\mathrm{Re} w - 1)} \Bigg\{ \Big\| \int_{\mathbb{R}^2} \frac{\mathbbm{1}_{t > |x - y|}}{\left(t^2 - |x - y|^2\right)^{\mathrm{Re} w}} |g(y)| \, dy \Big\|_{L^p\left(\mathbb{R}^2\right)} \nonumber \\
& \quad + \Big\| \int_{\mathbb{R}^2} \frac{\mathbbm{1}_{t > r_1 + r_2}}{(t^2 - \left(r_1 + r_2\right)^2)^{\mathrm{Re} w}} |g(y)| \, dy \Big\|_{L^p\left(\mathbb{R}^2\right)} \Bigg\} \nonumber \\
\lesssim& \|g\|_{L^p\left(\mathbb{R}^2\right)},
\end{align}
provided that
\begin{align}\label{estimate-xy}
\sup _{x \in \mathbb{R}^2} \int_{\mathbb{R}^2} \frac{\mathbbm{1}_{t>|x-y|}}{\left(t^2-|x-y|^2\right)^{\mathrm{Re} w}} \,dy \lesssim t^{2(1-\mathrm{Re} w)},
\end{align}
\begin{align}\label{estimate-r1r2}
\sup _{x \in \mathbb{R}^2} \int_{\mathbb{R}^2} \frac{\mathbbm{1}_{ t>r_1+r_2}}{(t^2-\left(r_1+r_2\right)^2)^{\mathrm{Re} w}}\,dy \lesssim t^{2(1-\mathrm{Re} w)}.
\end{align}

Indeed, we observe that $\frac{1}{2}<\mathrm{Re} w<1$, so this yields \eqref{shurtest} due to the Young inequality and by symmetry. Next, we proceed the proof of \eqref{estimate-xy} and \eqref{estimate-r1r2} briefly.

Estimate of \eqref{estimate-xy}:
$$
\int_{\mathbb{R}^2} \frac{\mathbbm{1}_{t>|x-y|}}{\left(t^2-|x-y|^2\right)^{\mathrm{Re} w}}\, d y \lesssim t^{-\mathrm{Re} w} \int_{|y|<t} \frac{1}{(t-|y|)^{\mathrm{Re} w}} \,dy \lesssim t^{2(1-\mathrm{Re} w)} .
$$

Estimate of \eqref{estimate-r1r2}:
$$
\int_{\mathbb{R}^2} \frac{\mathbbm{1}_{t>r_1+r_2}}{\left(t^2-(r_1+r_2)^2\right)^{\mathrm{Re} w}}\,dy \lesssim t^{-\mathrm{Re} w} \int_0^t \frac{r}{(t-r)^{\mathrm{Re} w}} \,d r \lesssim t^{2(1-\mathrm{Re} w)} .
$$

Thus, we show the conclusion of Theorem \ref{operator-family}.
\section{The proof of Theorem \ref{Lpestimate}}\label{proof-theorem}

In this section, we will focus on the proof of Theorem \ref{Lpestimate}. Inspired by the work of \cite{L}, the main idea is primarily based on the asymptotic properties of Bessel functions, which we will detail in the Appendix. Additionally, the proof utilizes some results regarding spectral multipliers, as stated in Lemma \ref{lem:multiplier}, within the context of magnetic fields.

\begin{figure}[htbp]
\begin{center}
\scalebox{0.9}[0.9]{
\begin{tikzpicture}
\draw (0,0) rectangle (6.75,6.75);
\draw[->]  (0,0) -- (0,7.5);
\draw[->]  (0,0) -- (7.5,0);
\draw (6.75,0) node[below] {$\ 1$};
\draw (0,0) node[below, left] {$0$};
\draw (7.5,0) node[above] {$1/p$};
\draw (0,7.5) node[right] {$l$};
\draw (3.35,0) circle (2pt) node[left] {$A$};    
\draw (6.75,1.7)  circle (2pt) node[right] {$B$}; 
\draw (0,1.7)  circle (2pt) node[left] {$C$}; 
\draw (6.75,3.35)  circle (2pt) node[right] {$\!D$}; 
\draw (0,3.35)  circle (2pt) node[left] {$\!E$}; 
\draw[dashed] (3.35,0) -- (0,1.7) ; 
\draw[dashed] (3.35,0) -- (6.75,1.7) ; 
\draw[dashed] (0,1.7) -- (6.75,1.7) ; 
\draw[dashed] (0,3.35) -- (6.75,3.35) ; 
\draw (3.35,0) node[below] {$\frac{1}{2}$};
\draw (5.8,0.8) node[below] {$l=\frac{1}{2}(\frac{1}{p}-\frac{1}{2})$};
\draw (1,0.8) node[below] {$l=\frac{1}{2}(\frac{1}{2}-\frac{1}{p})$};
\filldraw[fill=gray!50][dashed](0,1.7)--(6.75,1.7)--(3.35,0)--cycle; 
\filldraw[fill=blue!50](6.75,3.35)--(6.75,1.7)--(0,1.7)--(0,3.35);
\filldraw[fill=red!50](6.75,3.35)--(6.75,6.75)--(0,6.75)--(0,3.35);
\filldraw[fill=black](3.15,4.5) node[above right] {$\Omega_{3}$};
\filldraw[fill=black](3.15,2.25)  node[above right] {$\Omega_{2}$};
\filldraw[fill=black](3.15,0.7) node[above right] {$\Omega_{1}$};
\end{tikzpicture}
}
\end{center}
\caption{Here $A=(\frac{1}{2},0)$, $B=(1,\frac{1}{4})$, $C=(0,\frac{1}{4})$, $D=(1,\frac{1}{2})$, $E=(0,\frac{1}{2})$, respectively. The line $AB: l=\frac{1}{2}(\frac{1}{p}-\frac{1}{2})$. The line $AC: l=\frac{1}{2}(\frac{1}{2}-\frac{1}{p})$.}  \label{fig:1}
\end{figure}
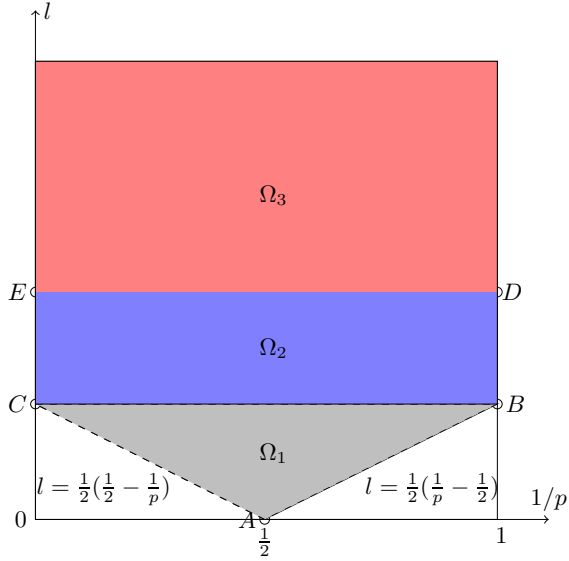

To prove Theorem \ref{Lpestimate}, it suffices to demonstrate that for all $\frac{1}{2}>\ell>\frac{1}{4}$ (corresponding to the region $\Omega_{2}$) and $1< p<+\infty$, there exists a constant $C(p, \ell)>0$ such that for all $f \in L^p\left(\mathbb{R}^2\right)$ and $t>0$
\begin{equation}\label{estimate-aim}
\left\|\left(1+\mathcal{L}_{\mathbf{A}}\right)^{-\ell} e^{i t \sqrt{\mathcal{L}_{\mathbf{A}}}} f\right\|_{L^p(\R^2)} \leq C(p, \ell)(1+t)^{2 \ell}\|f\|_{L^p(\R^2)} .
\end{equation}
This conclusion is illustrated in Figure \ref{fig:1}, where the gray, blue, red regions are labeled as $\Omega_{1}$, $\Omega_{2}$ and $\Omega_{3}$ respectively.
In fact, once we deduce that the region $\Omega_{2}$ is fulfilled, the other two regions are immediately available.

When $\ell \geq \frac{1}{2}$, it is evident that the operator $\left(1+\mathcal{L}_{\mathbf{A}}\right)^{-\left(\ell-\frac{3}{8}\right)}$ is bounded in $L^p\left(\mathbb{R}^2\right)$ for $1< p<+\infty$, with the norm that is independent of $p$. Combining with \eqref{estimate-aim}, this yields Theorem \ref{Lpestimate} for region $\Omega_{3}$.

For $\ell \leq \frac{1}{4}$, the point $A$ is trivial. Hence, in order to carry out Theorem \ref{Lpestimate} in region $\Omega_{1}$, we employ Stein's interpolation theorem (see \cite{stein}) to interpolate between the point $A$ and region $\Omega_{2}$.
\begin{remark}
As a matter of fact, the partition at $l\geq\frac{1}{2}$ is not strictly necessary, and $\frac{1}{2}$ can be substituted with any number larger than $\frac{1}{4}$. Accordingly, the boundedness of \eqref{estimate-aim} can also be altered to $(1+t)^{\frac{1}{2}+\varepsilon}$ for all $\varepsilon>0$.
\end{remark}

Finally, we focus on the proof of \eqref{estimate-aim} when $\frac{1}{2}>\ell>\frac{1}{4}$. Firstly, for all $s>0$, we define
$$
\begin{aligned}
&m(\ell, s)=\left(1+s^2\right)^{-\ell} e^{i s},\\
& m_{\ell}(s)=\psi(s) s^{-2 \ell} e^{i s}, \\
& M_{\ell}(s)=m(\ell, s)-m_{\ell}(s),
\end{aligned}
$$
in which $\psi \in C^{\infty}\left(\mathbb{R}_{+}\right)$ satisfies the following conditions:
\[
\psi(s) =
\begin{cases}
0, & \text{if } s \leq 1, \\
1, & \text{if } s \geq 2.
\end{cases}
\]
Next, we will show how the proof of \eqref{estimate-aim} can be reduced to establishing that for all $\frac{1}{2}>\ell>\frac{1}{4}$ and $1<p<+\infty$, there exists a constant $C(\ell, p)>0$ such that
\begin{equation}\label{estimate-mlt}
\|m(\ell, t \sqrt{\mathcal{L}_{\mathbf{A}}})\|_{p \rightarrow p} \leq C(\ell, p), \quad \forall t>0 .
\end{equation}
This result will serve as a fundamental step in our overall analysis.

Indeed, from \eqref{estimate:spectraltLA} in Lemma \ref{lem:multiplier}, we attain
$$
\begin{aligned}
\Big\|\left(1+\mathcal{L}_{\mathbf{A}}\right)^{-\ell} e^{i t \sqrt{\mathcal{L}_{\mathbf{A}}}}\Big\|_{p \rightarrow p} & \leq\Big\|\frac{\left(1+t^2 \mathcal{L}_{\mathbf{A}}\right)^{\ell}}{\left(1+\mathcal{L}_{\mathbf{A}}\right)^{\ell}}\Big\|_{p \rightarrow p}\times\Big\|\left(1+t^2 \mathcal{L}_{\mathbf{A}}\right)^{-\ell} e^{i t \sqrt{\mathcal{L}_{\mathbf{A}}}}\Big\|_{p \rightarrow p} \\
& \leq C \sum_{j=0}^2 \sup _{s \geq 0}\Big|s^j \frac{d^j}{d s^j} \frac{\left(1+s^2\right)^{\ell}}{\left(1+s^2/t^2\right)^{\ell}}\Big| \times \Big\|\left(1+t^2 \mathcal{L}_{\mathbf{A}}\right)^{-\ell} e^{i t \sqrt{\mathcal{L}_{\mathbf{A}}}}\Big\|_{p \rightarrow p}\\
&\lesssim (1+t)^{2\ell}\Big\|\left(1+t^2 \mathcal{L}_{\mathbf{A}}\right)^{-\ell} e^{i t \sqrt{\mathcal{L}_{\mathbf{A}}}}\Big\|_{p \rightarrow p}.
\end{aligned}
$$
Consequently, we can conclude that the truth of statement \eqref{estimate-mlt} guarantees the validity of Theorem \ref{Lpestimate}.
Then, it remains to show the following proposition.
\begin{proposition}\label{estimate:ml}
Let $\frac{1}{2}>\ell>\frac{1}{4}$ and $1<p<+\infty$. Then, there exists a constant $C(p, \ell)>0$ such that
\begin{equation}\label{estimate-ml}
\|m_{\ell}(t \sqrt{\mathcal{L}_{\mathbf{A}}}) f\|_{L^p(\R^2)} \leq C(p, \ell)\|f\|_{L^p(\R^2)}, \quad \forall f \in L^p\left(\mathbb{R}^2\right), t>0,
\end{equation}
and
\begin{equation}\label{estimate:Mls1}
\|M_{\ell}(t \sqrt{\mathcal{L}_{\mathbf{A}}}) f\|_{L^p(\R^2)} \leq C(\ell)\|f\|_{L^p(\R^2)}, \quad \forall f \in L^p\left(\mathbb{R}^2\right), t>0.
\end{equation}
\end{proposition}
\begin{proof}
The proof of Proposition \ref{estimate:ml} is analogous to the proof of Theorem 1.2 in \cite{L}. For the sake of clarity and fluency, we provide a brief proof.
Firstly, based on the spectral multiplier theorem (see Lemma \ref{lem:multiplier}), \eqref{estimate:Mls1} can be obtained directly.
We give the proof of \eqref{estimate:Mls1} shortly. By making full use of the Lagrange's mean value theorem, we observe that for all $j \in \mathbb{N}$, there exists a constant $C(j, \ell)>0$ such that
$$
\Big|\frac{d^j}{d s^j} M_{\ell}(s)\Big| \leq C(j, \ell)(1+s)^{-2 \ell-2}, \quad \forall s \geq 0 .
$$
So, by \eqref{estimate:spectraltLA}, we get
\begin{equation}\label{estimate:Mls}
\|M_{\ell}(t \sqrt{\mathcal{L}_{\mathbf{A}}})\|_{p \rightarrow p} \leq C(\ell), \quad \forall 1< p<+\infty, \ t>0 .
\end{equation}
Next, to complete the proof of \eqref{estimate-ml}, we consider
\begin{equation}\label{estimate:Besselcos}
\psi(s) s^{-\nu} \cos (s-\frac{\pi}{2} \nu)=\sqrt{\frac{\pi}{2}} s^{-(\nu-\frac{1}{2})} J_{\nu-\frac{1}{2}}(s)+W_{\nu}(s), \quad\forall s \geq 1,
\end{equation}
where $W_{\nu}(s)$ satisfies the estimate given by (refer to \eqref{rep:Besselcosrho} in the Appendix for more details)
\begin{equation}\label{estimate:Bessel}
|\frac{d^j}{d s^j} W_{\nu}(s)| \leq C(\nu, j)(1+s)^{-\mathrm{Re} \nu-1}, \quad\forall s \geq 0, \ j \in \mathbb{N} .
\end{equation}
We choose $\nu=2 \ell$ and $\nu=2 \ell+2 i$ respectively, in the context of \eqref{estimate:Besselcos}, leading to the following equalities:
\begin{align}\label{equality:cos}
\psi(s) s^{-2 \ell} \cos (s-\pi \ell) =\sqrt{\frac{\pi}{2}} s^{-(2 \ell-\frac{1}{2})} J_{2 \ell-\frac{1}{2}}(s)+W_{2 \ell}(s),
\end{align}
\begin{align}\label{equality:cosi}
\psi(s) s^{-2 \ell} \cos (s-\pi(\ell+i))  =\sqrt{\frac{\pi}{2}} s^{2i} s^{-(2 \ell+2 i-\frac{1}{2})} J_{2 \ell+2i-\frac{1}{2}}(s)+s^{2 i} W_{2 \ell+2i}(s) .
\end{align}

Note that
$$
\frac{e^{i \pi \ell}}{\sinh \pi}\left[\cos (s-\pi(\ell+i))-e^{-\pi} \cos (s-\pi \ell)\right]=e^{i s}.
$$
As a consequence, from \eqref{equality:cos} and \eqref{equality:cosi}, we can deduce
$$
\psi(s) s^{-2 \ell} e^{i s}=F(\ell, s)+N(\ell, s), \quad s>0,
$$
in which
\begin{equation}
F(\ell, s)=\frac{\sqrt{\frac{\pi}{2}}}{\sinh \pi}\left[e^{i \pi \ell} s^{2i} s^{-(2 \ell+2i -\frac{1}{2})} J_{2 \ell+2i-\frac{1}{2}}(s)-e^{\pi(i \ell-1)} s^{-(2 \ell-\frac{1}{2})} J_{2 \ell-\frac{1}{2}}(s)\right], \quad s>0,
\end{equation}
and
\begin{equation}
N(\ell, s)=\frac{1}{\sinh \pi}\left[e^{i \pi \ell} s^{2i} W_{2 \ell+2i}(s)-e^{\pi(i \ell-1)} W_{2 \ell}(s)\right], \quad s>0 .
\end{equation}
Let $w=(\frac32-2\ell)-2i$, then $\Re(w)=\frac32-2\ell:=\epsilon\in (1/2,1)$
for given $\frac{1}{4}<\ell<\frac{1}{2}$. Hence the assumption regarding $w$ in Theorem \ref{operator-family} is satisfied. Thus, we can apply Theorem \ref{operator-family} and \eqref{spectral-LA} in Lemma \ref{lem:multiplier} to obtain
$$
\|F(\ell, t \sqrt{\mathcal{L}_{\mathbf{A}}})\|_{p \rightarrow p} \leq C(p, \ell),\quad \forall 1<p<+\infty, t>0.
$$
We apply the Lemma \ref{lem:multiplier} again to derive the $L^p\rightarrow L^p$-boundedness of the operator $N(\ell, t \sqrt{\mathcal{L}_{\mathbf{A}}})$.
These results substantiate the assertions made in Proposition \ref{estimate:ml}.
\end{proof}

\section{The proof of Theorem \ref{sinLA}}\label{proof-sinLA}
The goal of this section is to show the following fixed time estimate for the solution of wave equation, i.e.
\begin{equation}
\left\|\frac{\sin \left(t \sqrt{\mathcal{L}_{\mathbf{A}}}\right)}{\sqrt{\mathcal{L}_{\mathbf{A}}}} f\right\|_{L^p\left(\mathbb{R}^2\right)} \leq C(p, t)\|f\|_{L^p\left(\mathbb{R}^2\right)}, \quad \forall p\geq1.
\end{equation}
Indeed, this can be proved by using the point-wise estimates of the Schwartz kernel of $\frac{\sin \left(t \sqrt{\mathcal{L}_{\mathbf{A}}}\right)}{\sqrt{\mathcal{L}_{\mathbf{A}}}}$, which was proved in \cite{FZZ}. We record it here for convenience.

\begin{lemma}\cite[Proposition 4.1]{FZZ}\label{sinkernelLA}
Let $K(t, x, y)$ be the Schwartz kernel of the operator $\frac{\sin \left(t \sqrt{\mathcal{L}_{\mathbf{A}}}\right)}{\sqrt{\mathcal{L}_{\mathbf{A}}}}$. Suppose $x=r_1\left(\cos \theta_1, \sin \theta_1\right)$ and $y=r_2\left(\cos \theta_2, \sin \theta_2\right)$, $\alpha$ to be in \eqref{alphadefine} and define
$$
\gamma=\frac{r_1^2+r_2^2-t^2}{2 r_1 r_2}=\frac{\left(r_1+r_2\right)^2-t^2}{2 r_1 r_2}-1=\frac{\left(r_1-r_2\right)^2-t^2}{2 r_1 r_2}+1
$$
and
$$
\beta_1=\cos ^{-1}\left(\frac{r_1^2+r_2^2-t^2}{2 r_1 r_2}\right), \quad \beta_2=\cosh ^{-1}\left(\frac{t^2-r_1^2-r_2^2}{2 r_1 r_2}\right).
$$
Then when $t \geq 0$, the kernel can be written as a ``geometric" term $G_w\left(t, r_1, \theta_1, r_2, \theta_2\right)$ and a ``diffractive" term $D_w\left(t, r_1, \theta_1, r_2, \theta_2\right)$
$$
K(t, x, y)=K\left(t, r_1, \theta_1, r_2, \theta_2\right)=G_w\left(t, r_1, \theta_1, r_2, \theta_2\right)+D_w\left(t, r_1, \theta_1, r_2, \theta_2\right),
$$
where
\begin{align}\label{Gwkernel}
& G_w\left(t, r_1, \theta_1, r_2, \theta_2\right) \nonumber\\
=&\frac{1}{2 \pi}\left(t^2-\left(r_1^2+r_2^2-2 r_1 r_2 \cos \left(\theta_1-\theta_2\right)\right)\right)^{-\frac{1}{2}} e^{i \int_{\theta_2}^{\theta_1} A\left(\theta^{\prime}\right) d \theta^{\prime}} \nonumber\\
& \quad\times\left\{\mathbbm{1}_{\left(\left|r_1-r_2\right|, r_1+r_2\right)}(t)\left[\mathbbm{1}_{\left[0, \beta_1\right]}\left(\left|\theta_1-\theta_2\right|\right)+e^{-2 \alpha \pi i} \mathbbm{1}_{\left[2 \pi-\beta_1, 2 \pi\right]}\left(\left|\theta_1-\theta_2\right|\right)\right]\right. \nonumber\\
&\left.\qquad+\mathbbm{1}_{\left(r_1+r_2, \infty\right)}(t)\left[\mathbbm{1}_{[0, \pi]}\left(\left|\theta_1-\theta_2\right|\right)+e^{-2 \alpha \pi i} \mathbbm{1}_{(\pi, 2 \pi]}\left(\left|\theta_1-\theta_2\right|\right)\right]\right\},
\end{align}
and
\begin{align}\label{Dwkernel}
& D_w\left(t, r_1, \theta_1, r_2, \theta_2\right)=\frac{\mathbbm{1}_{\left(r_1+r_2, \infty\right)}(t)}{\pi} e^{-i\left(\alpha\left(\theta_1-\theta_2\right)-\int_{\theta_2}^{\theta_1} A\left(\theta^{\prime}\right) d \theta^{\prime}\right)} \nonumber\\
& \quad \times \int_0^{\beta_2}\left(t^2-r_1^2-r_2^2-2 r_1 r_2 \cosh s\right)^{-\frac{1}{2}}\left(\sin (|\alpha| \pi) e^{-|\alpha| s}\right. \nonumber\\
& \left.+\sin (\alpha \pi) \frac{\left(e^{-s}-\cos \left(\theta_1-\theta_2+\pi\right)\right) \sinh (\alpha s)+i \sin \left(\theta_1-\theta_2+\pi\right) \cosh (\alpha s)}{\cosh (s)-\cos \left(\theta_1-\theta_2+\pi\right)}\right) \,d s.
\end{align}
When $t \leq 0$, the similar conclusion holds for \eqref{Gwkernel} and \eqref{Dwkernel} with replacing $t$ by $-t$.
\end{lemma}

\begin{lemma}[Pointwise estimate]\cite[Proposition 4.2]{FZZ}\label{estimate-point}
Let $G_w\left(t, r_1, \theta_1, r_2, \theta_2\right)$ be in \eqref{Gwkernel} and $D_w\left(t, r_1, \theta_1, r_2, \theta_2\right)$ be in \eqref{Dwkernel}. Then, in the polar coordinates $x=r_{1}(\cos\theta_{1},\sin\theta_{1})$, $y=r_{2}(\cos\theta_{2},\sin\theta_{2})$, the following estimates hold:
\begin{equation}
|G_w\left(t, r_1, \theta_1, r_2, \theta_2\right)|\lesssim \frac{1}{\sqrt{t^2-|x-y|^2}},\quad t^2>|x-y|^2,
\end{equation}
and
\begin{equation}
|D_w\left(t, r_1, \theta_1, r_2, \theta_2\right)|\lesssim \frac{1}{\sqrt{t^2-(r_{1}+r_{2})^2}},\quad t^2>(r_{1}+r_{2})^2.
\end{equation}
\end{lemma}
As a consequence, similar to the proof of Theorem \ref{operator-family}, we prove Theorem \ref{sinLA}
directly via combining Lemma \ref{sinkernelLA} and Lemma \ref{estimate-point}. In fact, to prove \eqref{est:sinLA}, it suffices to
show
\begin{equation}\label{estimate-sinxy}
\sup _{x \in \mathbb{R}^2} \int_{\mathbb{R}^2} \frac{\mathbbm{1}_{t>|x-y|}}{\left(t^2-|x-y|^2\right)^{\frac{1}{2}}} \,dy \lesssim t,\quad t>0,
\end{equation}
and
\begin{equation}\label{estimate-sinr1r2}
\sup _{x \in \mathbb{R}^2} \int_{\mathbb{R}^2} \frac{\mathbbm{1}_{ t>r_1+r_2}}{(t^2-\left(r_1+r_2\right)^2)^{\frac{1}{2}}}\,dy \lesssim t,\quad t>0.
\end{equation}
Therefore, a simple integral calculation,  as proving \eqref{estimate-xy} and \eqref{estimate-r1r2},  shows that \eqref{estimate-sinxy} and \eqref{estimate-sinr1r2} hold.

\appendix
\section{The Bessel function}
\addcontentsline{toc}{section}{Bessel's function}
In this section, we present the asymptotic behavior of the Bessel function $J_\kappa(s)$ as the variable $s$ approaches $+\infty$ with the fixed index $\kappa$, which we will use in the proof of Theorem \ref{Lpestimate}.

Let $H_\kappa^{1}(s)$ and $H_\kappa^{2}(s)$ denote the Hankel functions respectively. According to the definition of \cite[\S 3.61]{W}, we have
\begin{equation}\label{def:Bessel}
J_\kappa(s)=\frac{H_\kappa^{1}(s)+H_\kappa^{2}(s)}{2} .
\end{equation}

For $\mathrm{Re}(\kappa+\frac{1}{2})>0$ and $s>0$, the Hankel functions can be expressed as
$$
\begin{aligned}
& H_\kappa^{1}(s)=\sqrt{\frac{2}{\pi s}} \frac{e^{i\left(s-\frac{1}{2} \pi \kappa-\frac{\pi}{4}\right)}}{\Gamma\left(\kappa+\frac{1}{2}\right)} \int_0^{\infty} e^{-u} u^{\kappa-\frac{1}{2}}\left(1+i \frac{u}{2 s}\right)^{\kappa-\frac{1}{2}} \,d u, \\
& H_\kappa^{2}(s)=\sqrt{\frac{2}{\pi s}} \frac{e^{-i\left(s-\frac{1}{2} \pi \kappa-\frac{\pi}{4}\right)}}{\Gamma\left(\kappa+\frac{1}{2}\right)} \int_0^{\infty} e^{-u} u^{\kappa-\frac{1}{2}}\left(1-i \frac{u}{2 s}\right)^{\kappa-\frac{1}{2}} \,d u.
\end{aligned}
$$
These integral representations can be found in \cite[\S 6.12]{W} by setting $\beta=0$. The first representation is also available in \cite[\S 8.421, p.915]{GR}, and the seconde one can be derived from the first and the relation $H_\kappa^{2}(s)=-e^{i \pi \kappa} H_\kappa^{1}(-s)$, as shown in the 8th equality of \cite[\S 8.476, p.927]{W}).

We define
$$
\begin{aligned}
\left(\frac{1}{2}-\kappa\right)_0=1,\left(\frac{1}{2}-\kappa\right)_j=\prod_{i=0}^{j-1}\left(\frac{1}{2}-\kappa+i\right), \quad j \in \mathbb{N}^*.
\end{aligned}
$$
Additionally, let
$$
\begin{aligned}
R_\kappa(k, i w)=\int_0^1(1-h)^{k-1}(1+i w h)^{\kappa-\frac{1}{2}-k} d h, \quad w \in \mathbb{R}, \ k \in \mathbb{N}^* .
\end{aligned}
$$
Then, for $k \in \mathbb{N}^*$, we can express $(1+i w)^{\kappa-\frac{1}{2}}$ as follows (see also \cite[\S 7.2]{W})
$$
(1+i w)^{\kappa-\frac{1}{2}}=\sum_{j=0}^{k-1} \frac{\left(\frac{1}{2}-\kappa\right)_j}{j!}\left(\frac{w}{i}\right)^j+\frac{\left(\frac{1}{2}-\kappa\right)_k}{(k-1)!}\left(\frac{w}{i}\right)^k R_\kappa(k, i w), \quad w \in \mathbb{R} .
$$

Therefore, we can obtain
\begin{align}\label{def:Hankle}
H_\kappa^{1}(s)= & \sqrt{\frac{2}{\pi s}} \frac{e^{i\left(s-\frac{1}{2} \pi \kappa-\frac{\pi}{4}\right)}}{\Gamma\left(\kappa+\frac{1}{2}\right)}\left[\sum_{j=0}^{k-1} \frac{\left(\frac{1}{2}-\kappa\right)_j \Gamma\left(\kappa+j+\frac{1}{2}\right)}{j!(i 2 s)^j}\right. \nonumber\\
& \left.+\frac{\left(\frac{1}{2}-\kappa\right)_k}{(k-1)!(i 2 s)^k} \int_0^{\infty} e^{-u} u^{\kappa+k-\frac{1}{2}} R_\kappa\left(k, i \frac{u}{2 s}\right) \,d u\right],
\end{align}
\begin{align}\label{def:Hankles}
H_\kappa^{2}(s)= & \sqrt{\frac{2}{\pi s}} \frac{e^{-i\left(s-\frac{1}{2} \pi \kappa-\frac{\pi}{4}\right)}}{\Gamma\left(\kappa+\frac{1}{2}\right)}\left[\sum_{j=0}^{k-1} \frac{\left(\frac{1}{2}-\kappa\right)_j \Gamma\left(\kappa+j+\frac{1}{2}\right)}{j!(-i 2 s)^j}\right. \\\nonumber
& \left.+\frac{\left(\frac{1}{2}-\kappa\right)_k}{(k-1)!(-i 2 s)^k} \int_0^{\infty} e^{-u} u^{\kappa+k-\frac{1}{2}} R_\kappa\left(k,-i \frac{u}{2 s}\right) \,d u\right] .
\end{align}
Consequently, for $\kappa \in \mathbb{C}$ with $\mathrm{Re} \kappa>-\frac{1}{2}$, by \eqref{def:Bessel}-\eqref{def:Hankles}, we can compute
\begin{equation}\label{Bessel-cosrep}
s^{-\kappa} J_\kappa(s)=\sqrt{\frac{2}{\pi}} s^{-\left(\kappa+\frac{1}{2}\right)} \cos (s-\frac{\pi}{2}(\kappa+\frac{1}{2}))+R_\kappa(s), \quad s \geq \frac{1}{2},
\end{equation}
where $R_\kappa(s)$ satisfies the following estimate
\begin{equation}\label{Bessel-error}
\Big|\frac{d^j}{d s^j} R_\kappa(s)\Big| \leq C(k, j) s^{-\mathrm{Re} \kappa-\frac{3}{2}}, \quad \forall s \geq 1,\ j \in \mathbb{N} .
\end{equation}

Let $\psi \in C^{\infty}\left(\mathbb{R}_{+}\right)$ be a smooth function that satisfies the following conditions
$$
0 \leq \psi \leq 1, \psi(h)=0 \text { if } h \leq 1 \text {, and } \psi(h)=1 \text { if } h \geq 2 .
$$

Then, for $\mathrm{Re} \nu>0$, by utilizing \eqref{Bessel-cosrep} and \eqref{Bessel-error}, along with the fact that the function $z^{-\kappa} J_\kappa(z)(z \in \mathbb{C})$ is analytic, we can derive the following expression
\begin{equation}\label{rep:Besselcosrho}
\psi(s) s^{-\nu} \cos (s-\frac{\pi}{2} \nu)=\sqrt{\frac{\pi}{2}} s^{-\left(\nu-\frac{1}{2}\right)} J_{\nu-\frac{1}{2}}(s)+W_{\nu}(s),\quad \forall s \geq 0,
\end{equation}
where $W_{\nu}(s)$ satisfies the following estimate:
$$
\Big|\frac{d^j}{d s^j} W_{\nu}(s)\Big| \leq C(\nu, j)(1+s)^{-\mathrm{Re} \nu-1}, \quad \forall s \geq 0,\ j \in \mathbb{N} .
$$

\begin{center}

\end{center}
 \end{document}